\RequirePackage{fix-cm}
\RequirePackage{amsmath}
\documentclass[smallcondensed,final]{svjour3}     
\smartqed  

\usepackage{amsfonts}
\usepackage{graphicx}
\usepackage{transparent}
\usepackage{romannum}
\usepackage[gen]{eurosym}
\usepackage{enumerate}
\usepackage{array}
\usepackage{empheq}
\usepackage{color}
\usepackage{tcolorbox}
\usepackage{siunitx}
\usepackage[toc,page]{appendix}

\graphicspath{{images/}}

\newcommand{\pder}[2]{\frac{\partial #1}{\partial #2}}

\DeclareMathOperator{\tr}{tr}
\begin{document}

\title{Extension of B-spline Material Point Method for unstructured triangular grids using Powell-Sabin splines}

\titlerunning{Powell-Sabin spline based Material Point Method}        

\author{Pascal de Koster \and
		Roel Tielen \and 
		Elizaveta Wobbes \and
		Matthias M\"oller        
}


\institute{P. de Koster \at
              Department of Applied Mathematics, Section
of Numerical Analysis, Delft University of
Technology, Van Mourik Broekmanweg 6,
2628 XE Delft, The Netherlands \\
              \email{p.b.j.dekoster@tudelft.nl}           
}

\date{Received: date / Accepted: date}

\maketitle

\begin{abstract} The Material Point Method (MPM) is a numerical technique that combines a fixed Eulerian background grid and Lagrangian point masses to simulate materials which undergo large deformations. 
Within the original MPM, discontinuous gradients of the piecewise-linear basis functions lead to so-called `grid-crossing errors' when particles cross element boundaries. Previous research has shown that B-spline MPM (BSMPM) is a viable alternative not only to MPM, but also to more advanced versions of the method that are designed to reduce the grid-crossing errors. In contrast to many other MPM-related methods, BSMPM has been used exclusively on structured rectangular domains, considerably limiting its range of applicability. In this paper, we present an extension of BSMPM to unstructured triangulations.  The proposed approach combines MPM with $C^1$-continuous high-order Powell-Sabin (PS) spline basis functions. Numerical results demonstrate the potential of these basis functions within MPM in terms of grid-crossing-error elimination and higher-order convergence.
\end{abstract}

\section{Introduction}\label{sec:Introduction}
The Material Point Method (MPM) has proven to be
successful in solving complex
engineering problems that involve large deformations, multi-phase interactions and
history-dependent material behaviour. 
Over the years, MPM has been applied to  a wide range of applications, including modelling of failure phenomena in single- and multi-phase media \cite{andersen2010land,Sulsky2004,Zabala2011}, crack growth  \cite{Nairn2003,guo2017simulation} and snow and ice dynamics  \cite{Stomakhin2013,sulsky2007ice,Gaume2018}.
Within MPM, a continuum is discretised by defining a set of Lagrangian particles, called material points, which store all relevant material properties. 
An Eulerian background grid is adopted on which the equations of motion are solved in every time step. The solution on the background grid is then used to update all material-point properties such as displacement, velocity, and stress. 
In this way, MPM incorporates both Eulerian and Lagrangian descriptions. Similarly to other combined Eulerian-Lagrangian techniques, MPM attempts to avoid the numerical difficulties arising from nonlinear convective terms associated with an Eulerian problem formulation, while at the same time preventing grid distortion, typically encountered within mesh-based Lagrangian formulations (e.g.,~\cite{Donea2004,sulsky1994particle}).

Classically, MPM uses piecewise-linear Lagrange basis functions, also known as `tent' functions. However, the gradients of these basis functions are discontinuous at element boundaries. This leads to so-called grid-crossing errors~\cite{bardenhagen2004generalized} when material points cross this discontinuity. Grid-crossing errors can significantly influence the quality of the numerical solution and may eventually lead to a lack of convergence (e.g.,~\cite{steffen2008analysis}). Different methods have been developed to mitigate the effect of grid-crossings. For example, the Generalised Interpolation Material Point (GIMP) \cite{bardenhagen2004generalized} and Convected Particle Domain Interpolation (CPDI)~\cite{Sadeghirad2011} methods eliminate grid-crossing inaccuracies by introducing an alternative particle representation. The GIMP method represents material points by particle-characteristic functions and reduces to standard MPM, when the Dirac delta function centered at the material-point position is selected as the characteristic function. For multivariate cases, a number of versions of the GIMP method are available such as the uniform GIMP (uGIMP) and contiguous-particle GIMP (cpGIMP) methods. The CPDI method extends GIMP in order to accurately capture shear distortion. Much research has been performed to further improve the accuracy of the CPDI approach and increase its range of applicability~\cite{Sadeghirad2013,Nguyen2017,Leavy2019}. On the other hand, the Dual Domain Material Point (DDMP) method~\cite{zhang2011material} preserves the original point-mass representation of the material points, but adjusts the gradients of the basis functions to avoid grid-crossing errors.
The DDMP method replaces the gradients of the piecewise-linear Lagrange basis functions in standard MPM by smoother ones. The DDMP method with sub-points~\cite{Dhakal2016} proposes an alternative manner for numerical integration within the DDMP algorithm.

The B-spline Material Point Method (BSMPM)~\cite{steffen2008analysis,Steffen2010} solves the problem of grid-crossing errors completely by replacing piecewise-linear Lagrange basis functions with higher-order B-spline basis functions. 
B-spline and piecewise-linear Lagrange basis functions possess many common properties. For instance, they both satisfy the partition of unity, are non-negative and have a compact support. The main advantage of higher-order B-spline basis functions over piecewise-linear Lagrange basis function is, however, that they have at least $C^0$-continuous gradients which preclude grid-crossing errors from the outset. Moreover, spline basis functions are known to provide higher accuracy per degree of freedom as compared to $C^0$-finite elements~\cite{Hughes2014}. On structured rectangular grids, adopting B-spline basis functions within MPM not only eliminates grid-crossing errors but also yields higher-order spatial convergence~\cite{steffen2008analysis,Steffen2010,tielen2017high,Gan2018mpm,Bazilevs2017}.  Previous research also demonstrates that BSMPM is a viable alternative to the GIMP, CPDI and  DDMP methods~\cite{steffen2008examination,Motlagh2017,Gan2018mpm,wobbestaylor}. While the CPDI and DDMP methods can be used on unstructured grids~\cite{zhang2011material,Nguyen2017,Leavy2019}, to the best of our knowledge, BSMPM for unstructured grids does not yet exist. This implies that its applicability to real-world problems is limited compared to the CPDI and DDMP methods. 

In this paper, we propose an extension of BSMPM to unstructured triangulations to combine the benefits of B-splines with the geometric flexibility of triangular grids. The proposed method employs quadratic Powell-Sabin (PS) splines~\cite{powell1977piecewise}. These splines are
piecewise higher-order polynomials defined on a particular refinement of any given triangulation and are typically used in computer aided geometric design and approximation theory~\cite{speleers2012isogeometric,dierckx1992algorithms,Manni2007,Sablonniere1987}. 
PS splines are $C^1$-continuous and hence overcome the grid-crossing issue within MPM by design. We would like to remark that although this paper focuses on PS splines, other options such as refinable $C^1$ splines~\cite{Nguyen2016} can be used to extend MPM to unstructured triangular grids. The proposed PS-MPM approach is analysed based on several benchmark problems. 

The paper is organised as follows. In Section~\ref{sec:MPM}, the governing equations are provided together with the MPM solution strategy. In Section~\ref{sec:PS-splines}, the construction of PS-spline basis functions and their application within MPM is described. In Section~\ref{sec:ResultsAndDiscussion}, the obtained numerical results are presented. In Section~\ref{sec:Conclusion}, conclusions and recommendation are given.



\section{Material Point Method}
\label{sec:MPM}

We will summarise the MPM as introduced by Sulsky et al. \cite{sulsky1994particle} to keep this work self-contained. First, the governing equations are presented, afterwards the MPM strategy to solve these equations is introduced. 
\subsection{Governing equations}\label{subsec:GoverningEquations}
The deformation of a continuum is modelled using the conservation of linear momentum and a material model. It should be noted that MPM can be implemented with a variety of material models that either use the rate of deformation (i.e., the symmetric part of the velocity gradient) or the deformation gradient. However, for this study it is sufficient to consider the simple linear-elastic and neo-Hookean models that are based on the deformation gradient. Using the Einstein summation convention, the system of equations in a Lagrangian frame of reference for each direction $x_k$ is given by
\begin{align}
    \pder{u_k}{t} &= v_k, \label{eq:Displacement}\\
    \rho\pder{v_k}{t} &= \pder{\sigma_{kl}}{x_l}+\rho g_k,     	\label{eq:MomentumConservation}\\
   D_{kl} &= \delta_{kl} + \pder{u_k}{x^0_l},\\
   \sigma_{kl} &= \begin{cases} 
    \lambda \delta_{kl} \tr\left(E_{kl}\right) + 2\mu E_{kl} & \text{for linear-elastic materials,}\\
   \frac{\lambda \ln(J)}{J}\delta_{kl}+\frac{\mu}{J}\left(D_{km}D_{lm} - \delta_{kl} \right)  & \text{for neo-Hookean materials.} 
   \end{cases}
   \label{eq:Stress}
\end{align}
where $u_k$ is the displacement, $t$ is the time, $v_k$ is the velocity, $\rho$ is the density, $g_k$ is the body force, $D_{kl}$ is the deformation gradient, $\delta_{kl}$ is the Kronecker delta, $\sigma_{kl}$ is the stress tensor, $x^0$ is the position in the reference configuration, $\lambda$ is the Lam\'{e} constant, $E_{kl}$ is the strain tensor defined as $\frac{1}{2}\left(D_{kl} + D_{lk} \right)- \delta_{kl}$, $J$ is the determinant of the deformation gradient and $\mu$ is the shear modulus. Equation~\ref{eq:MomentumConservation} and Equation~\ref{eq:Stress} describe the conservation of linear momentum and the material model, respectively.

Initial conditions are required for the displacement, velocity and stress tensor. The boundary of the domain $\Omega$ can be divided into a part with a Dirichlet boundary condition for the displacement and a part with a Neumann boundary condition for the traction:
\begin{alignat}{2}
    u_k(\boldsymbol{x},t) &=U_k(\boldsymbol{x},t)&\quad\text{on }\partial\Omega_u,\\
    \sigma_{kl}(\boldsymbol{x},t)\, n_l &= \tau_k(\boldsymbol{x},t)&\quad\text{on }\partial\Omega_\tau,
\end{alignat}
where $\boldsymbol{x} = [\:x_1\:\: x_2\:]^T$ and $\boldsymbol{n}$ is the unit vector normal to the boundary $\partial\Omega$ and pointing outwards. We remark that the domain $\Omega$, its boundary $\partial\Omega = \partial\Omega_u\cup \partial\Omega_\tau$ and the normal unit vector $\boldsymbol{n}$ may all depend on time.

For solving the equations of motion in the Material Point Method, the conservation of linear momentum in Equation~\ref{eq:MomentumConservation} is needed in its weak form. For the weak form, Equation~\ref{eq:MomentumConservation} is first multiplied by a continuous test function $\phi$ that vanishes on $\partial \Omega_u$, and is subsequently integrated over the domain $\Omega$:   
\begin{equation}\label{eq:WeakMomentumEquationTemp}
	\int_\Omega\phi\rho a_k \,\mathrm{d}\Omega=\int_{\Omega} \phi \pder{\sigma_{kl}}{x_l}\, \mathrm{d}\Omega + \int_\Omega\phi\rho g_k \,\mathrm{d}\Omega,
\end{equation}
in which $a_k=\pder{v_k}{t}$ is the acceleration. After applying integration by parts, the Gauss integration theorem and splitting the boundary into Dirichlet and Neumann part, the weak formulation becomes
\begin{equation}\label{eq:WeakMomentumEquation}
	\int_\Omega\phi\rho a_k \,\mathrm{d}\Omega=
	\int_{\partial\Omega_\tau} \phi \tau_k\, \mathrm{d}S \\
	-\int_\Omega \pder{\phi}{x_l} \sigma_{kl}\, \mathrm{d}\Omega 
	+ \int_\Omega\phi\rho g_k \,\mathrm{d}\Omega,
\end{equation}
whereby the contribution on $\partial\Omega_u$ equals zero, as the test function vanishes on this part of the boundary. 

\subsection{Discretised equations}\label{subsec:DiscretisedEquations}

Equation~\ref{eq:WeakMomentumEquation} can be solved using a finite element approach by defining $n_{bf}$ basis functions $\phi_i$ $(i = 1, \dots, n_{bf})$. The acceleration field $a_k$ is then discretised as a linear combination of these basis functions:
\begin{equation}\label{eq:AccelerationProjection}
	a_k(\boldsymbol{x},t) = \sum_{j=1}^{n_{bf}} \hat{a}_{k,j}\phi_j(\boldsymbol{x}),
\end{equation}
in which $\hat{a}_{k,j}$ is the time-dependent $j^{th}$ acceleration coefficient corresponding to basis function $\phi_j$. Substituting Equation~\ref{eq:AccelerationProjection} into Equation~\ref{eq:WeakMomentumEquation} and expanding the test function in terms of the basis functions $\phi_i$ leads to
\begin{equation}
	\int_\Omega\phi_i\rho \left(\sum_{j=1}^{n_{bf}} \hat{a}_{k,j}\phi_j\right) \,\mathrm{d}\Omega= 
	\int_{\partial\Omega_\tau} \phi_i \tau_k\, \mathrm{d}S \\
	- \int_\Omega \pder{\phi_i}{x_l} \sigma_{kl}\, \mathrm{d}\Omega 
	+ \int_\Omega\phi_i\rho g_k \,\mathrm{d}\Omega, 
\end{equation}
which holds for $i=1,\dots, n_{bf}$. By exchanging summation and integration, this can be rewritten in matrix-vector form as follows:
\begin{equation}\label{eq:SystemOfEqs0}
	\sum_{j=1}^{n_{bf}}\underbrace{\left(\int_\Omega\phi_i\rho\phi_j\,\mathrm{d}\Omega\right)}_{M_{ij}}
	\hat{a}_{k,j}=
	\underbrace{\int_{\partial\Omega_\tau} \phi_i \tau_k\, \mathrm{d}S}_{F_{k,i}^{trac}}\\
	- \underbrace{\int_\Omega \pder{\phi_i}{x_l} \sigma_{kl}\, \mathrm{d}\Omega}_{F_{k,i}^{int}}
	+ \underbrace{\int_\Omega\phi_i\rho g_k \,\mathrm{d}\Omega}_{F_{k,i}^{body}},
\end{equation}
\begin{equation}\label{eq:SystemOfEqs}
	\Rightarrow\boldsymbol{M}\boldsymbol{\hat{a}}_{k}=
	\boldsymbol{F}_k^{trac} 
	-\boldsymbol{F}_k^{int} 
	+\boldsymbol{F}_k^{body},
\end{equation}
where $\boldsymbol{M}$ denotes the mass matrix, $\boldsymbol{\hat{a}}_{k}$ the coefficient vector for the acceleration, while $\boldsymbol{F}_k^{trac}$, $\boldsymbol{F}_k^{int}$ and $\boldsymbol{F}_k^{body}$ denote respectively the traction, internal and body force vector in the $x_k$-direction. When density, stress and body force are known at time $t$, the coefficient vector $\boldsymbol{\Hat{a}}_k$ is found from Equation~\ref{eq:SystemOfEqs}.  The properties of the continuum at time $t+\Delta t$ can then be determined using the acceleration field at time $t$.

\subsection{Solution procedure}\label{subsec:MPM}
\begin{figure}[h tb]
	\centering
	\begin{minipage}[b]{0.3\textwidth}
		\includegraphics[width=\textwidth]{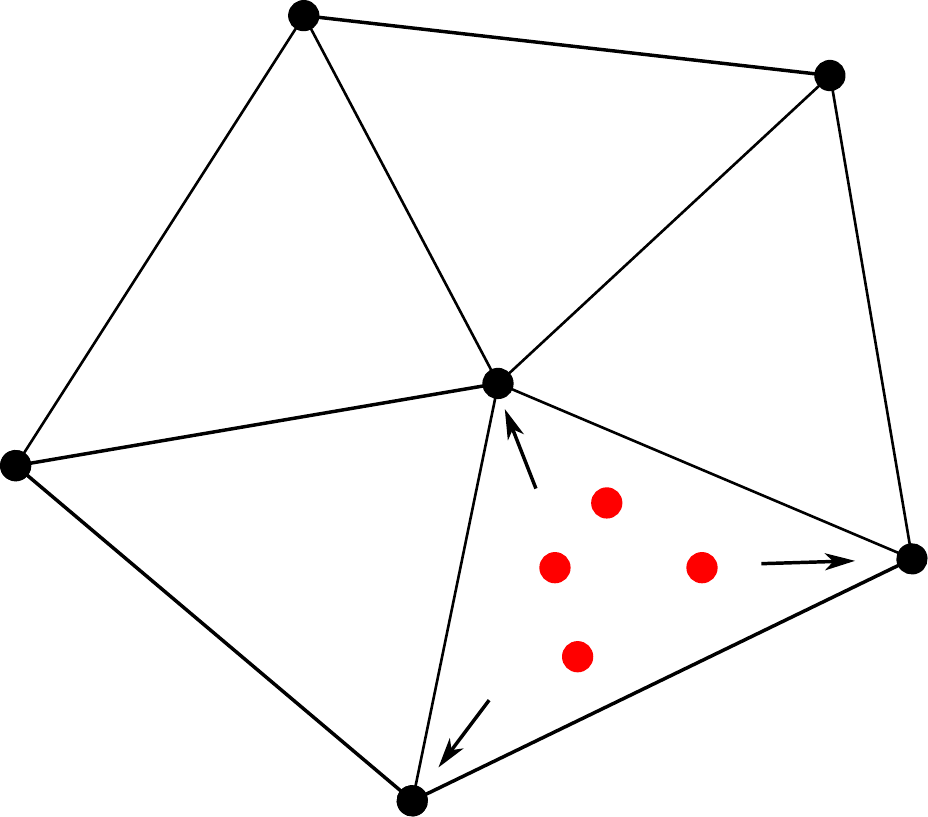}
	\end{minipage}
	\begin{minipage}[b]{0.3\textwidth}
		\includegraphics[width=\textwidth]{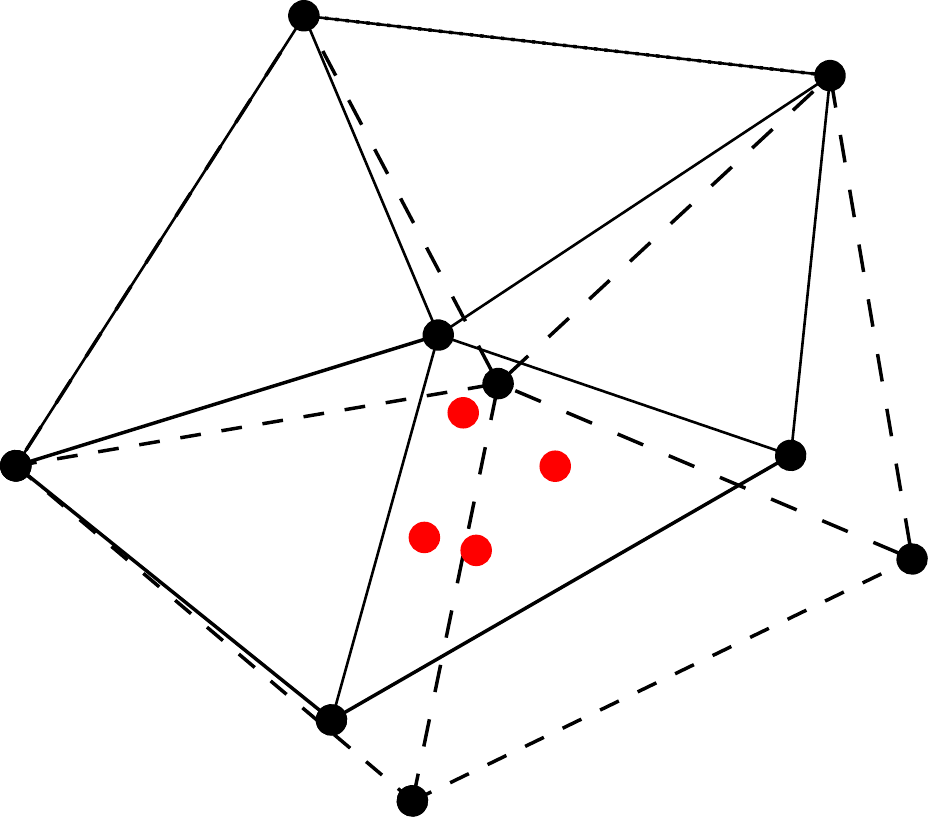}
	\end{minipage}
	\begin{minipage}[b]{0.3\textwidth}
		\includegraphics[width=\textwidth]{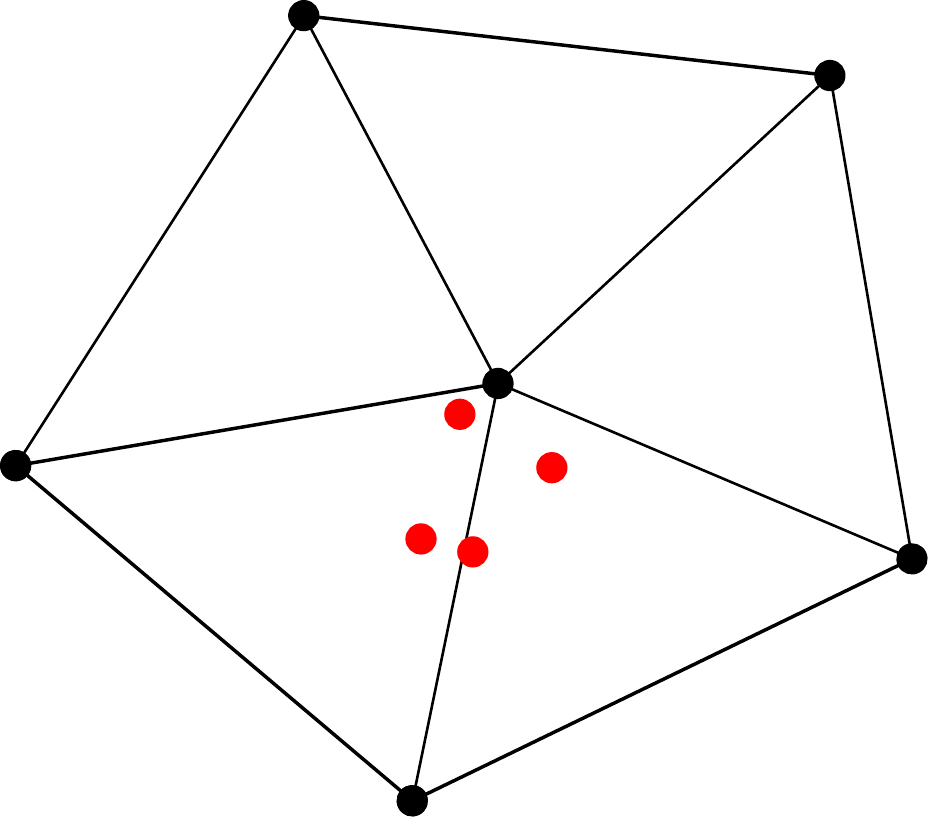}
	\end{minipage}
	\caption{One time step in the MPM procedure. Projection of particle properties onto the basis functions of the background grid (left). Update of particle properties based on the solution on the background grid (middle). Resetting the grid (right).}
	\label{fig:MPMGrid}
\end{figure}

Within MPM, the continuum is discretised by a set of particles that store all its physical properties. At each time step, the particle information is projected onto a background grid, on which the momentum equation is solved. Particle properties and positions are updated according to this solution, as illustrated in Figure~\ref{fig:MPMGrid}.

To solve Equation~\ref{eq:SystemOfEqs} at every time step, the integrals in Equation~\ref{eq:SystemOfEqs0} have to be evaluated. Since the material properties are only known at the particle positions, these positions are used as integration points for the numerical quadrature, and the particle volumes $V$ are adopted as integration weights. Assume that the total number of particles is equal to $n_p$. Here, subscript $p$ corresponds to particle properties. A superscript $t$ is assigned to particle properties that change over time. The mass matrix and force vectors are then defined by
\begin{align}
	M^t_{ij}=
	&\sum_{p=1}^{n_p} \underbrace{V_p^t \rho_p^t}_{m_p} \phi_i(\boldsymbol{x}_p^t) \phi_j(\boldsymbol{x}_p^t), \label{eq:MassMatrixAssembled}\\
	F^{trac,t}_{k,i} =
	& \sum_{p=1}^{n_{p}} \phi_i(\boldsymbol{x}_{p}^t)\, f^{trac,t}_{k,p},\label{eq:TracForceAssembled}\\
	F^{int,t}_{k,i} =
	&  \sum_{p=1}^{n_p} V_p^t  \pder{\phi_i}{x_m}(\boldsymbol{x}_p^t)\,\sigma_{mk,p}^t,\label{eq:IntForceAssembled}\\
	F^{body,t}_{k,i} =
	& \sum_{p=1}^{n_p} \underbrace{V_p^t \rho_p^t}_{m_p} \phi_i(\boldsymbol{x}_p^t) g_k, \label{eq:extForceAssembled}
\end{align}
in which $m_p$ denotes the particle mass, which is set to remain constant over time, guaranteeing conservation of the total mass. The coefficient vector $\boldsymbol{\hat{a}}_k^t$ at time $t$ for the $x_k$-direction is determined by solving 
\begin{equation}\label{eq:SystemOfEqs2}
	\boldsymbol{M}^{t}\boldsymbol{\hat{a}}_{k}^t=\left(\boldsymbol{F}_k^{trac,t} -\boldsymbol{F}_k^{int,t} +\boldsymbol{F}_k^{body,t}\right).
\end{equation}
Unless otherwise stated, Equation~\ref{eq:SystemOfEqs2} is solved adopting the consistent mass matrix. Next, the reconstructed acceleration field is used to update the particle velocity:
\begin{equation}\label{eq:VelocityUpdate}
	v_{k,p}^{t+\Delta t}=v_{k,p}^{t}+\Delta t \,  a_k^t(\boldsymbol{x}_p^t) =
	v_{k,p}^{t}+\Delta t\sum_{j=1}^{n_{bf}} \hat{a}_{k,j}^{t} \phi_j(\boldsymbol{x}_p^t).
\end{equation}

The semi-implicit Euler-Cromer scheme \cite{cromer1981stable} is adopted to update the remaining particle properties. First, the velocity field $v_k^{t+\Delta t}$ is discretised as a linear combination of the same basis functions, similar to the acceleration field (Eq.~\ref{eq:AccelerationProjection}):
\begin{equation}\label{eq:VelocityFieldProjection}
	v_k^{t+\Delta t}(\boldsymbol{x})= \sum_{j=1}^{n_{bf}} \hat{v}_{k,j}^{t+\Delta t}\phi_j(\boldsymbol{x}),
\end{equation}
in which $\boldsymbol{\hat{v}}_{k}^{t+\Delta t}$ is the velocity coefficient vector for the velocity field at time $t+\Delta t$ in the $x_k$-direction. The coefficient vector $\boldsymbol{\hat{v}}_{k}^{t+\Delta t}$ is then determined from a density weighted $L_2$-projection onto the basis functions $\phi_i$. This results in the following system of equations:
\begin{equation}\label{eq:SystemVelocityProjection}
	\boldsymbol{M}^t\boldsymbol{\hat{v}}^{t+\Delta t}_k=\boldsymbol{P}^{t+\Delta t}_k,
\end{equation}
where $\boldsymbol{M}^t$ is the same mass matrix as defined in Equation~\ref{eq:MassMatrixAssembled}. $\boldsymbol{P}_{k}^{t+\Delta t}$ denotes the momentum vector with the coefficients given by
\begin{equation}
	{P}^{t+\Delta t}_{k,i} = \sum_{p=1}^{n_p} \underbrace{V_p^t \rho_p^t}_{m_p} v_{k,p}^{t+\Delta t}\phi_i(\boldsymbol{x}_p).
\end{equation}

Particle properties are subsequently updated in correspondence with Equations~\ref{eq:Displacement}-\ref{eq:Stress}. First, the deformation gradient and its determinant are obtained:
\begin{align}   
   D_{kl,p}^{t+\Delta t}=& 
    \left(\delta_{kl}+\varepsilon_{kl,p}^{t+\Delta t}\Delta t\right)D_{kl,p}^t,\\
    J_p^{t+\Delta t}=& \det(\boldsymbol{D}_p^{t+\Delta t}).
\end{align}
Here, $\varepsilon_{kl,p}^{t+\Delta t}$ denotes the symmetric part of the velocity gradient:
\begin{equation}
    \varepsilon_{kl,p}^{t+\Delta t}=
    \frac{1}{2}\sum_{j=1}^{{n_{bf}}}\left(
 	\hat{v}_{k,j}^{t+\Delta t}
 	\pder{\phi_j}{x_l}(\boldsymbol{x}_p^t)+
 	\hat{v}_{l,j}^{t+\Delta t}
 	\pder{\phi_j}{x_k}(\boldsymbol{x}_p^t) \right).
\end{equation}
For linear-elastic materials, the particle stresses are computed as follows:
\begin{equation}
    E_{kl,p}^{t+\Delta t} = \frac{1}{2}\left(D_{kl,p}^{t+\Delta t}+D_{lk,p}^{t+\Delta t}\right) - \delta_{kl},
\end{equation}
\begin{equation} \label{eq:StressIncrement1}
\sigma_{kl,p}^{t+\Delta t}=
    \lambda \delta_{kl} \tr\left( E_{kl,p}^{t+\Delta t}\right) + 2\mu E_{kl,p}^{t+\Delta t}.
\end{equation}
For neo-Hookean materials, the stresses are obtained from
\begin{equation} \label{eq:StressIncrement2}
\sigma_{kl,p}^{t+\Delta t}=    
    \frac{\lambda \ln(J_p^{t+\Delta t})}{J_p^{t+\Delta t}}\delta_{kl} + \frac{\mu}{J_p^{t+\Delta t}}\left(D_{km,p}^{t+\Delta t}D_{lm,p}^{t+\Delta t}-\delta_{kl}\right).
\end{equation}

The determinant of the deformation gradient is used to update the volume of each particle. Based on this volume, the density is updated in such way that the mass of each particle remains constant:
\begin{align}
    V_p^{t+\Delta t}=& J_p^{t+\Delta t}V_p^0,\\
    \rho_p^{t+\Delta t}=& m_p/V_p^{t+\Delta t}.
\end{align}

Finally, particle positions and displacements are updated from the velocity field:
\begin{align}
	x_{k,p}^{t+\Delta t}=& x_{k,p}^{t} +\Delta t \sum_{i=1}^{{n_{bf}}} v_{k,i}^{t+\Delta t} \phi_i(\boldsymbol{x}_p^t), \label{eq:PositionUpdate}\\
	u_{k,p}^{t+\Delta t}=& u_{k,p}^{t} +\Delta t \sum_{i=1}^{{n_{bf}}} v_{k,i}^{t+\Delta t} \phi_i(\boldsymbol{x}_p^t). \label{eq:DisplacementUpdate}
\end{align}
The described MPM can be numerically implemented by performing the steps in Equations~\ref{eq:MassMatrixAssembled}-\ref{eq:DisplacementUpdate} in each time step in the shown order. Note that all steps can be applied with a variety of basis functions, without any essential difference in the algorithm. In this paper, we investigate the use of Powell-Sabin (PS) spline basis functions, which are described in the next section.


\section{Powell-Sabin spline basis functions}\label{sec:PS-splines}

\begin{figure}[h!tb]
	\center{
	\begin{minipage}[b]{0.48\textwidth}
	    \centering
	    \includegraphics[width=\textwidth]{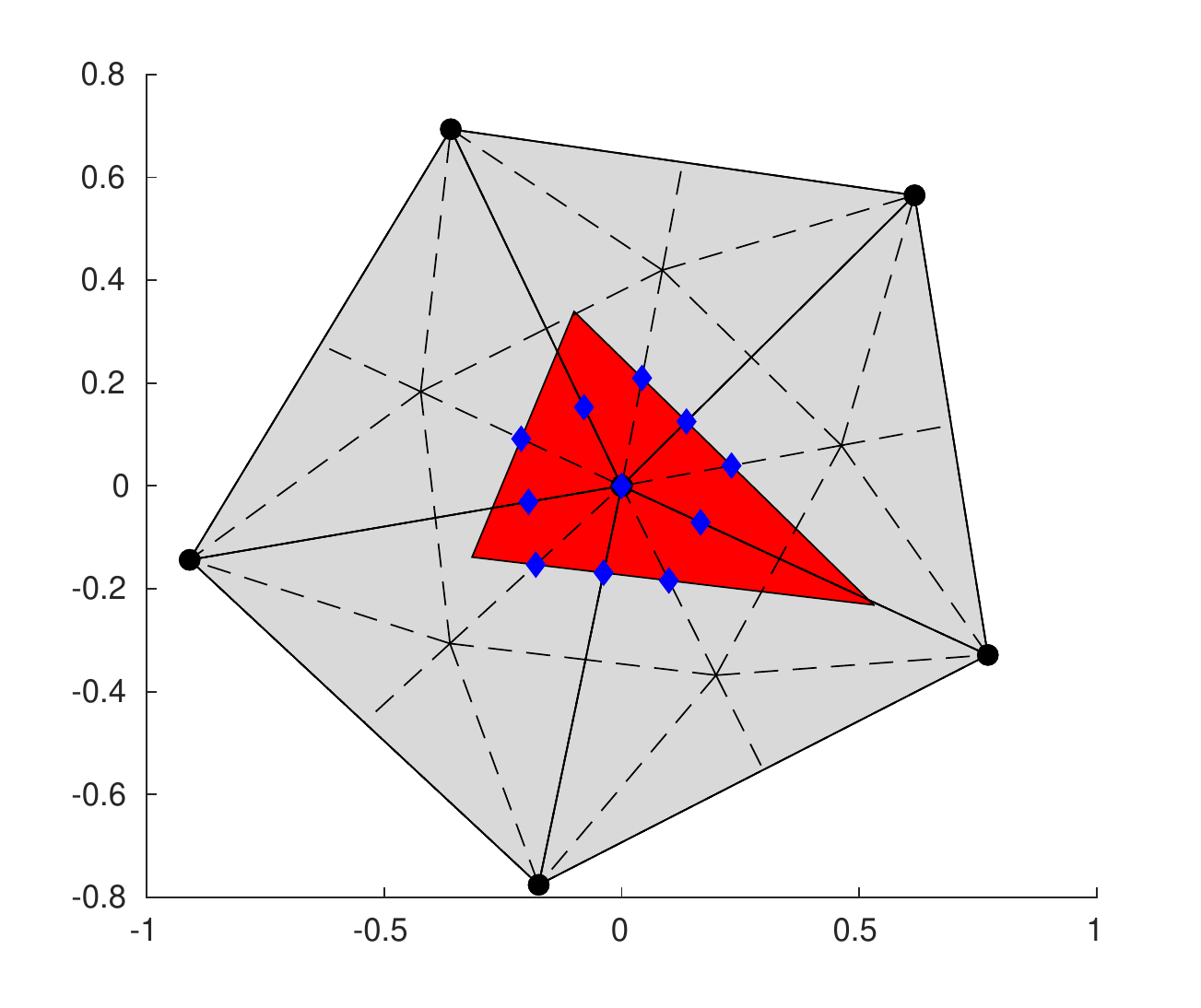}
	\end{minipage}
	\begin{minipage}[b]{0.48\textwidth}
	    \centering
	    \includegraphics[width=\textwidth]{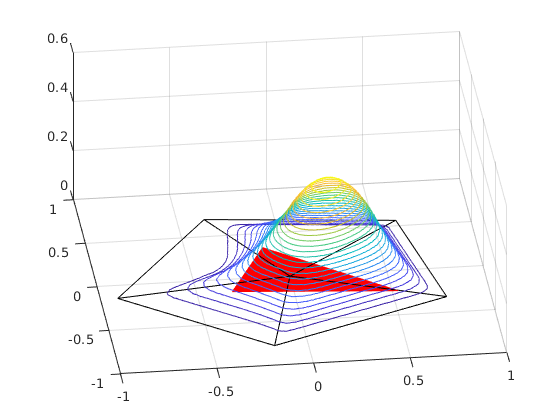}
	\end{minipage}
	}
	\center{
	\begin{minipage}[b]{0.48\textwidth}
	    \centering
	    \includegraphics[width=\textwidth]{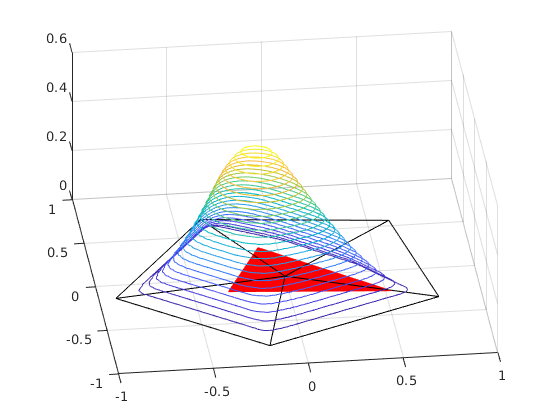}
	\end{minipage}
	\begin{minipage}[b]{0.48\textwidth}
		\centering
	    \includegraphics[width=\textwidth]{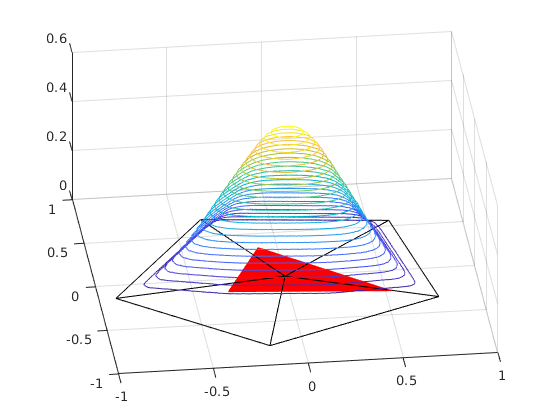}
	\end{minipage}
	}
	\caption{The PS-points for the middle vertex marked with diamonds and a control triangle (top left). The three PS-splines associated with this control triangle are shown in the other subfigures.}
    \label{fig:PS-spline}
\end{figure}

In this section, the essential properties of Powell-Sabin (PS) spline basis functions are shortly summarised, for the full construction process and the proofs of the properties we refer to the publication by Dierckx~\cite{dierckx1997calculating,dierckx1992algorithms}. The considered PS-splines are piecewise quadratic splines with global $C^1$-continuity. They are defined on arbitrary triangulations, have local support and possess the properties of non-negativity and partition of unity. Several examples of PS-splines are shown in Figure~\ref{fig:PS-spline}.

\begin{figure}[h!tb]
    \begin{minipage}[b]{0.48\textwidth}
	    \centering
		\def\svgwidth{1\columnwidth}
\begingroup%
  \makeatletter%
  \providecommand\color[2][]{%
    \errmessage{(Inkscape) Color is used for the text in Inkscape, but the package 'color.sty' is not loaded}%
    \renewcommand\color[2][]{}%
  }%
  \providecommand\transparent[1]{%
    \errmessage{(Inkscape) Transparency is used (non-zero) for the text in Inkscape, but the package 'transparent.sty' is not loaded}%
    \renewcommand\transparent[1]{}%
  }%
  \providecommand\rotatebox[2]{#2}%
  \newcommand*\fsize{\dimexpr\f@size pt\relax}%
  \newcommand*\lineheight[1]{\fontsize{\fsize}{#1\fsize}\selectfont}%
  \ifx\svgwidth\undefined%
    \setlength{\unitlength}{453.17224727bp}%
    \ifx\svgscale\undefined%
      \relax%
    \else%
      \setlength{\unitlength}{\unitlength * \real{\svgscale}}%
    \fi%
  \else%
    \setlength{\unitlength}{\svgwidth}%
  \fi%
  \global\let\svgwidth\undefined%
  \global\let\svgscale\undefined%
  \makeatother%
  \begin{picture}(1,0.94886537)%
    \lineheight{1}%
    \setlength\tabcolsep{0pt}%
    \put(0,0){\includegraphics[width=\unitlength,page=1]{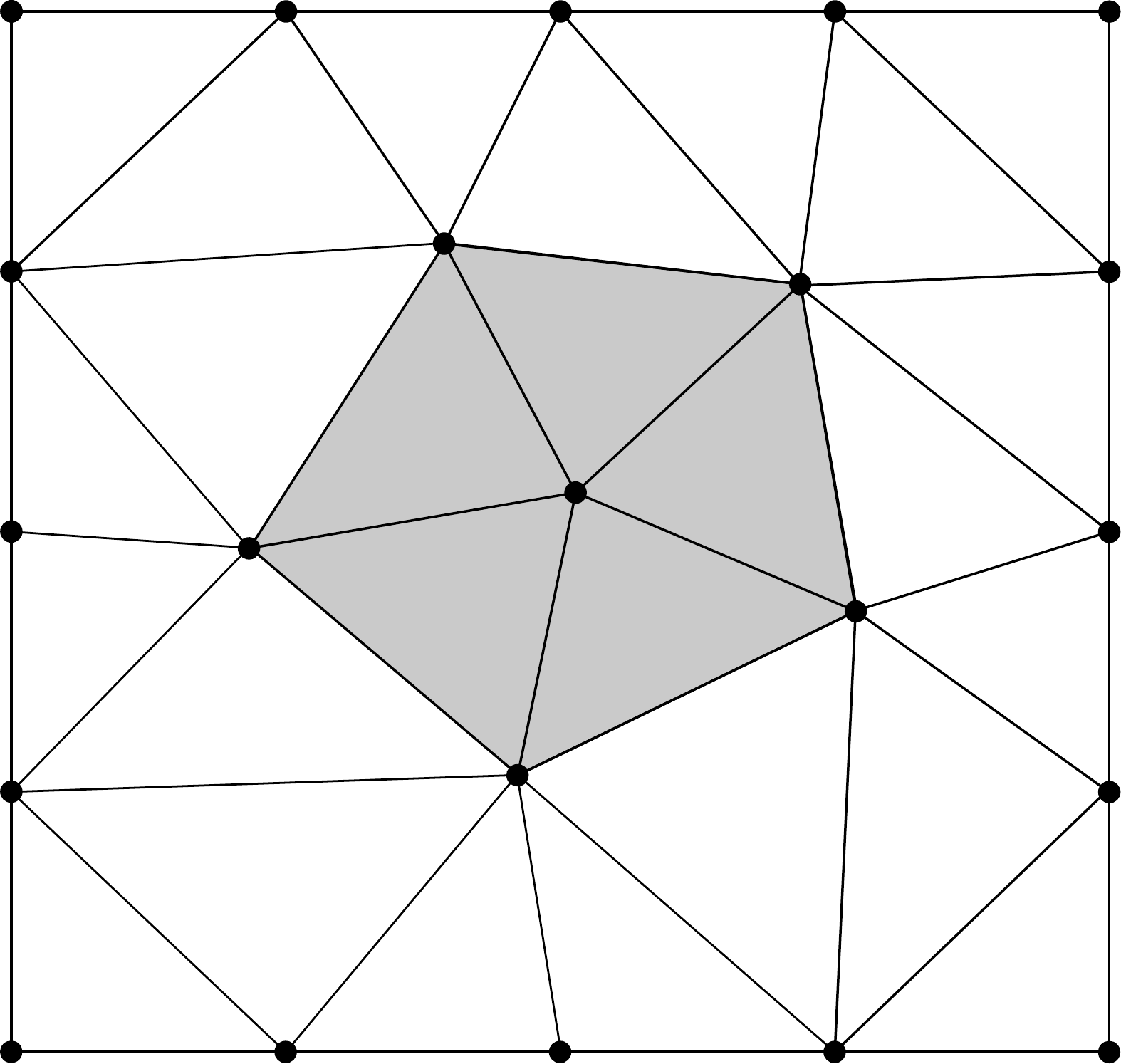}}%
    \put(0.34860943,0.5056924){\color[rgb]{0,0,0}\makebox(0,0)[lt]{\lineheight{0}\smash{\begin{tabular}[t]{l}$t_j$\end{tabular}}}}%
    \put(0.56984482,0.61412333){\color[rgb]{0,0,0}\makebox(0,0)[lt]{\lineheight{0}\smash{\begin{tabular}[t]{l}$t_m$\end{tabular}}}}%
    \put(0.51762047,0.44799658){\color[rgb]{0,0,0}\makebox(0,0)[lt]{\lineheight{0}\smash{\begin{tabular}[t]{l}$V_i$\end{tabular}}}}%
  \end{picture}%
\endgroup%

	\end{minipage}
	\hspace{0.02\textwidth}
	\begin{minipage}[b]{0.48\textwidth}
	    \centering
		\def\svgwidth{1\columnwidth}
		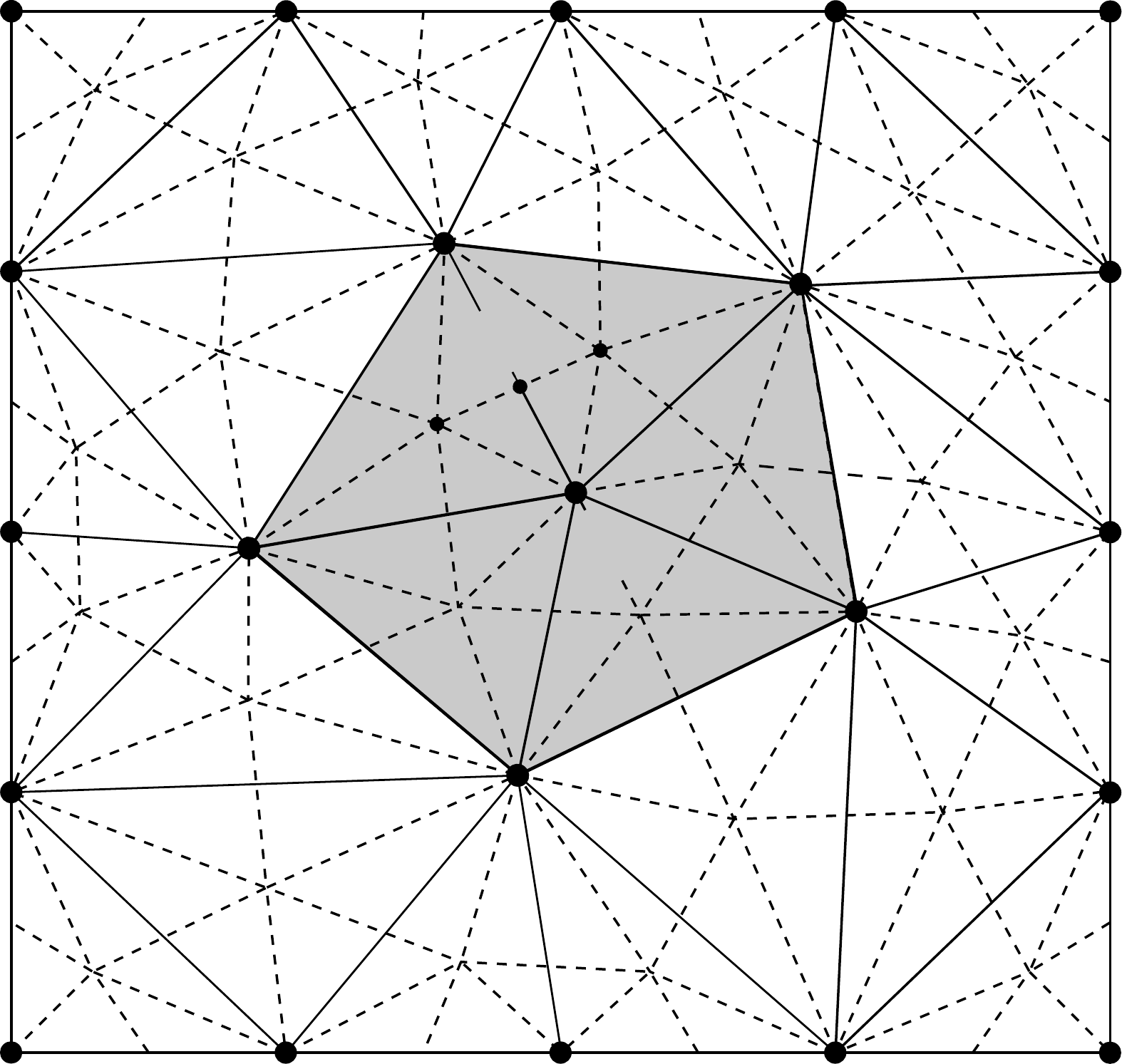
	\end{minipage}
	\caption{A triangular grid (left) and a possible PS-refinement (right). The molecule $\Omega_i$ around vertex $V_i$ is marked grey.}
	\label{fig:PSRefinement}
\end{figure}

To construct PS-splines on an arbitrary triangulation, a Powell-Sabin refinement is required, dividing each of the original main triangles into six sub-triangles as follows (see also Figure~\ref{fig:PSRefinement}).  
\begin{enumerate}
	\item For each triangle $t_j$, choose an interior point $Z_j$ (e.g., the incenter), such that if triangles $t_j$ and $t_m$ have a common edge, the line between $Z_j$ and $Z_m$ intersects the edge. The intersection point is called $Z_{jm}$.
	\item Connect each $Z_j$ to the vertices of its triangle $t_j$ with straight lines.
	\item Connect each $Z_j$ to all edge points $Z_{jm}$ with straight lines. In case $t_j$ is a boundary element, connect $Z_j$ to an arbitrary point on each boundary edge (e.g., the edge middle).
\end{enumerate}


Each PS-spline is associated with a main vertex $V_i$, and will be locally defined on its \textit{molecule} $\Omega_i$, which is defined as the area consisting of the main triangles around $V_i$, as shown in Figure~\ref{fig:PSRefinement}. Next, for each $V_i$, a set of PS-points is defined as the union of $V_i$ and the sub-triangle edge midpoints directly around it (see Figure~\ref{fig:PS-spline}). A control triangle for $V_i$ is then defined as a triangle that contains all its PS-points, and preferably has a small surface \cite{dierckx1997calculating}. The set of control triangles uniquely defines the set of PS-splines, in which three PS-splines are associated with each control triangle or main vertex. The construction process from the control triangles ensures the properties of non-negativity and partition of unity (the details are omitted here, but can be found in \cite{dierckx1997calculating,dierckx1992algorithms}). Figure~\ref{fig:PS-spline} shows a PS-refinement, a main vertex and its associated molecule, a possible control triangle and its three associated PS-spline basis functions. Note that for PS-splines, each main vertex is associated with three PS-spline basis functions, in contrast to a piecewise-linear Lagrange basis, for which each main vertex is associated with only one basis function.

\section{Numerical Results}\label{sec:ResultsAndDiscussion}

To validate the proposed PS-MPM, several benchmarks involving large deformations are considered. The first benchmark describes a thin vibrating bar, where the displacement is caused by an initial velocity. For this benchmark, PS-MPM on a structured triangular grid is compared with a reference solution. After that, a vibrating plate with known analytical solution is considered. The spatial convergence of PS-MPM on an unstructured grid is determined for this benchmark.
The third benchmark describes a soil-column under self-weight. For this benchmark, we propose the strategy of partial lumping to mitigate oscillations in the numerical solution that occur when the mass matrix is lumped. 

\subsection{Vibrating bar}\label{sec:VibratingBar}

\begin{figure}[tbh]
	\centering
	\begin{minipage}{1\linewidth}
		\centering
		\def\svgwidth{0.95\textwidth}
\begingroup%
  \makeatletter%
  \providecommand\color[2][]{%
    \errmessage{(Inkscape) Color is used for the text in Inkscape, but the package 'color.sty' is not loaded}%
    \renewcommand\color[2][]{}%
  }%
  \providecommand\transparent[1]{%
    \errmessage{(Inkscape) Transparency is used (non-zero) for the text in Inkscape, but the package 'transparent.sty' is not loaded}%
    \renewcommand\transparent[1]{}%
  }%
  \providecommand\rotatebox[2]{#2}%
  \newcommand*\fsize{\dimexpr\f@size pt\relax}%
  \newcommand*\lineheight[1]{\fontsize{\fsize}{#1\fsize}\selectfont}%
  \ifx\svgwidth\undefined%
    \setlength{\unitlength}{392.60365074bp}%
    \ifx\svgscale\undefined%
      \relax%
    \else%
      \setlength{\unitlength}{\unitlength * \real{\svgscale}}%
    \fi%
  \else%
    \setlength{\unitlength}{\svgwidth}%
  \fi%
  \global\let\svgwidth\undefined%
  \global\let\svgscale\undefined%
  \makeatother%
  \begin{picture}(1,0.31610697)%
    \lineheight{1}%
    \setlength\tabcolsep{0pt}%
    \put(0.46588644,0.28278301){\color[rgb]{0,0,0}\makebox(0,0)[lt]{\lineheight{0}\smash{\begin{tabular}[t]{l}$v_0$\end{tabular}}}}%
    \put(0,0){\includegraphics[width=\unitlength,page=1]{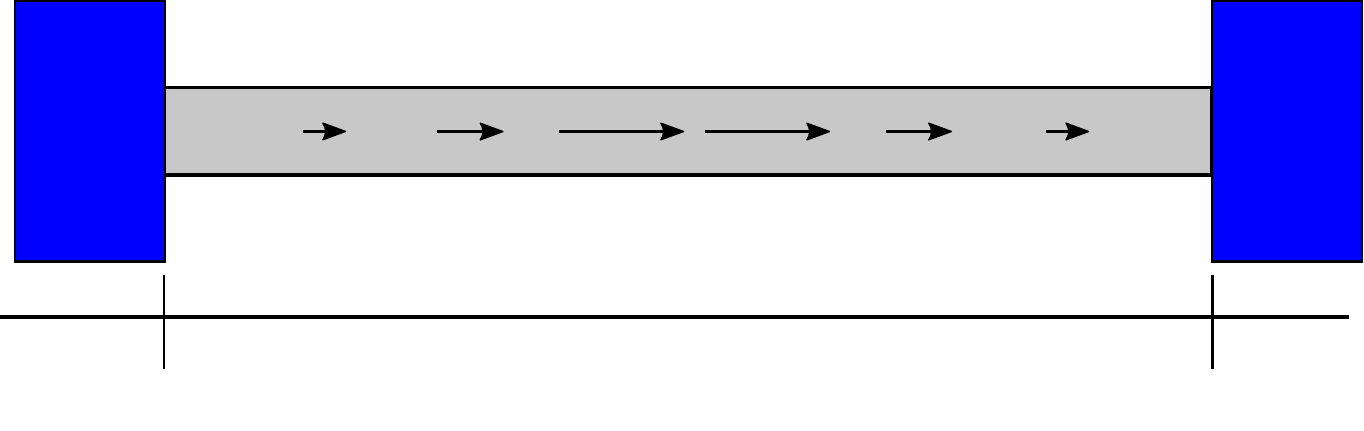}}%
    \put(0.10169317,0.00530716){\color[rgb]{0,0,0}\makebox(0,0)[lt]{\lineheight{0}\smash{\begin{tabular}[t]{l}$0$\end{tabular}}}}%
    \put(0.87204847,0.01050934){\color[rgb]{0,0,0}\makebox(0,0)[lt]{\lineheight{0}\smash{\begin{tabular}[t]{l}$L$\end{tabular}}}}%
  \end{picture}%
\endgroup%

	\end{minipage}
	\begin{minipage}{0.80\textwidth}
	    \hspace*{-3mm} \includegraphics[width=\textwidth]{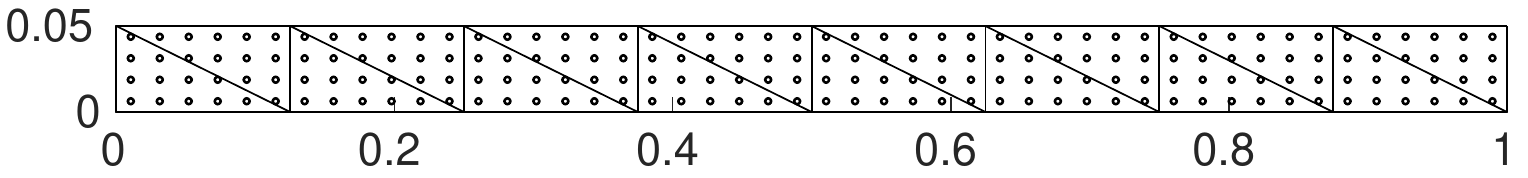}
	\end{minipage}
    \caption{The vibrating bar benchmark (top) and the initial particle configuration (bottom).}
    \label{fig:GridBsplineBar}
\end{figure}
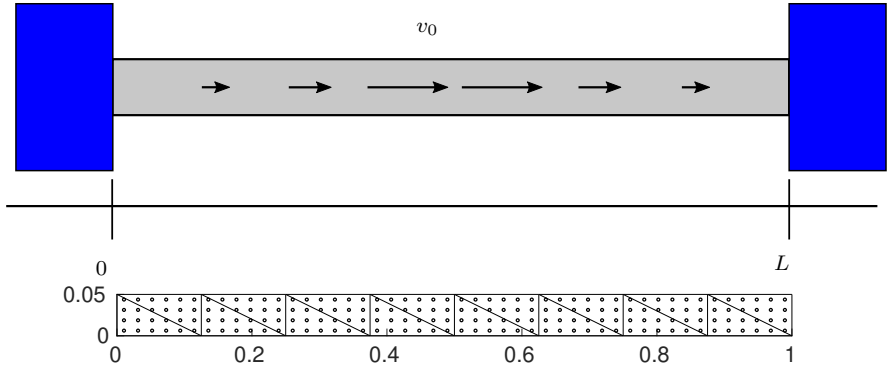

In this section, a thin linear-elastic vibrating bar is considered with both ends fixed. A UL-FEM \cite{ten2007large} solution on a very fine grid serves as reference. The grid used for PS-MPM and the initial particle positions are shown in Figure~\ref{fig:GridBsplineBar}. 

The bar is modelled with density $\rho=25\,\si{kg/m^3}$, Young's modulus $E=50\,\si{Pa}$, Poisson ratio $\nu = 0$, length $L=1\, \si{m}$, width $W=2\,\si{m}$ and initial maximum velocity $v_0=0.1\, \si{m/s}$. The chosen parameters result in a maximum normal strain in the $x$-direction of approximately $7\%$. At the left and right boundary, homogeneous Dirichlet boundary conditions are imposed for both $x$- and $y$-displacement, whereas at the top and bottom boundary, a homogeneous Dirichlet boundary condition is imposed only for the $y$-displacement, and a free-slip boundary condition for the $x$-direction.
The initial displacement equals zero, but an initial $x$-velocity profile is set to $v_x(x^0,y^0,t) = 0.1\sin(x^0\pi/L)$. The initial $y$-velocity is equal to zero. 

The time step size for the simulation is $\Delta t = 5\cdot 10^{-3}$ s. The corresponding Courant number is defined as $C = \tfrac{\Delta t}{h}\sqrt{E/\rho}$, in which $\sqrt{E/\rho}$ is the characteristic wave speed of the material and $h$ is the typical element length. Due to the ambiguity of $h$ for a PS-refinement on an unstructured triangular grid, we estimate the average edge length of the sub-triangles by $h\approx 0.025$~m. In this case, the Courant number is $C\approx 0.57< 1$, satisfying the CFL condition \cite{courant1928partiellen}. 

\begin{figure}[tbh]
	\centering
	\begin{minipage}{0.45\textwidth}
		\centering
		\includegraphics[width=\textwidth]{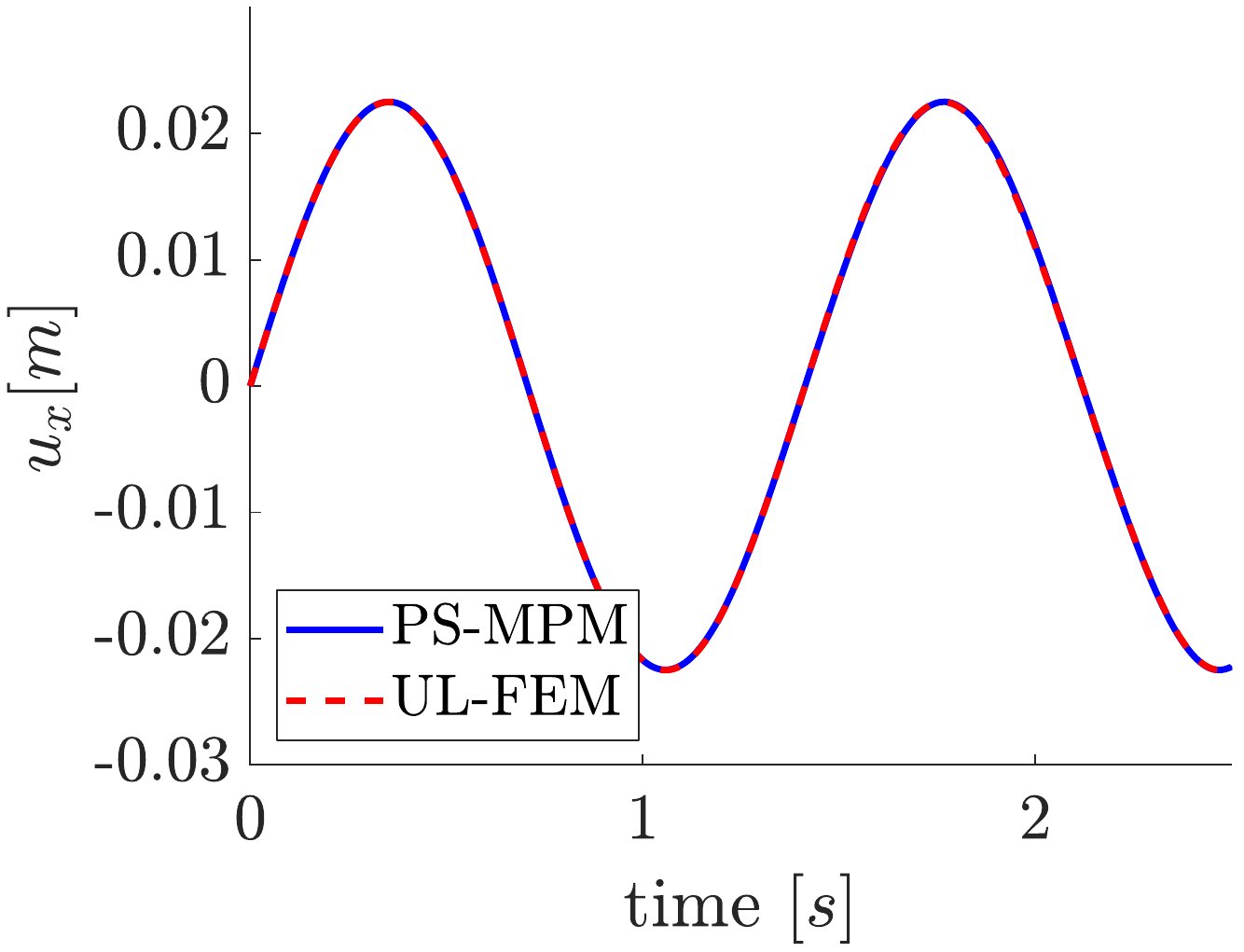}
	\end{minipage}
	\hspace{0.2cm}
	\begin{minipage}{0.45\textwidth}
	    \includegraphics[width=\textwidth]{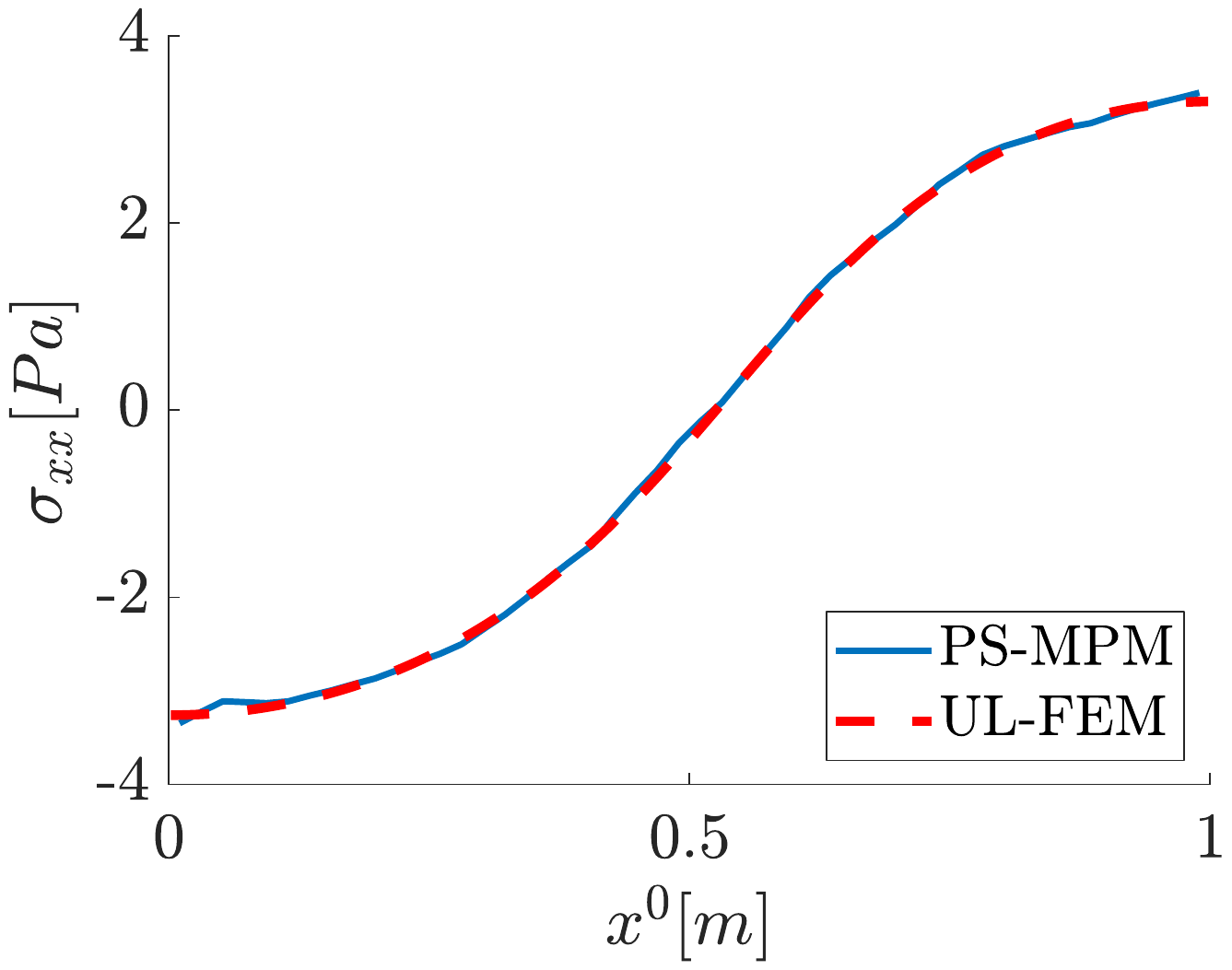}
    \end{minipage}
    \caption{Displacement of the particle initially positioned at $(x^0,y^0) \approx (0.5$~m,$0.025$~m$)$ over time (left) and the stress component $\sigma_{xx}$ along the bar at $t=2.5$ s (right).}
    \label{fig:VibBar}
\end{figure}

Figure~\ref{fig:VibBar} shows the displacement in the middle of the bar over time and the stress profile through the entire bar at the end of the simulation. Although relatively coarse PS-MPM grid and particle configurations are adopted, the method yields accurate results and shows a smooth stress profile.

\subsection{Vibrating plate undergoing axis-aligned displacement}

In this section, a neo-Hookean vibrating plate on the unit square ($L \times W$, with $L=W=1$ m) is considered. The plate undergoes axis-aligned displacement. An analytic solution for this problem is constructed using the method of manufactured solutions (MMS): the analytic solution is assumed \textit{a priori}, and the corresponding body force is calculated accordingly. This benchmark has been adopted from \cite{sadeghirad2011convected}. The analytical solution for the displacement in terms of the reference configuration is given by
\begin{align}
    u_x &= u_0\sin\left(2\pi x^0\right)\sin\left(\sqrt{E/\rho_0}\,\pi t\right), \label{eq:MMS_ux}\\
    u_y &= u_0\sin\left(2\pi y^0\right)\sin\left(\sqrt{E/\rho_0}\,\pi t+\pi\right)\label{eq:MMS_uy}.
\end{align}
Here, $\rho_0=10^3\,\si{kg/m^3}$, $u_0 =0.05\,\si{m}$ and $E=10^7\,\si{Pa}$. The corresponding body forces \cite{sadeghirad2011convected} are
\begin{align}
    g_x &= \pi^2 u_x\left(\frac{4\mu}{\rho_0}-\frac{E}{\rho_0}-4\frac{\lambda \left[\ln (D_{xx}D_{yy})-1\right]-\mu}{\rho_0 D_{xx}^2}\right),\\
    g_y &= \pi^2 u_y\left(\frac{4\mu}{\rho_0}-\frac{E}{\rho_0}-4\frac{\lambda \left[\ln (D_{xx}D_{yy})-1\right]-\mu}{\rho_0 D_{yy}^2}\right), 
\end{align}
in which the Lam\'e constant $\lambda$, the shear modulus $\mu$, and the normal components of the deformation gradient $D_{xx}$ and $D_{yy}$ are defined as
\begin{align}
    \lambda &= \frac{E\nu}{(1+\nu)(1-2\nu)},\quad\quad \mu = \frac{E}{2(1+\nu)},\\
    D_{xx} &= 1+2u_0\pi \cos(2\pi x^0)\sin\left(\sqrt{E/\rho_0}\,\pi t\right) \label{eq:MMS_Fxx},\\
    D_{yy} &= 1+2u_0\pi \cos(2\pi y^0)\sin\left(\sqrt{E/\rho_0}\,\pi t+\pi\right) \label{eq:MMS_Fyy}.
\end{align}
Here, $\nu = 0.3$ and all solutions are again given with respect to the reference configuration. 

\begin{figure}[htb]
\centering
\begin{minipage}[b]{0.49\textwidth}
	\includegraphics[width=\textwidth]{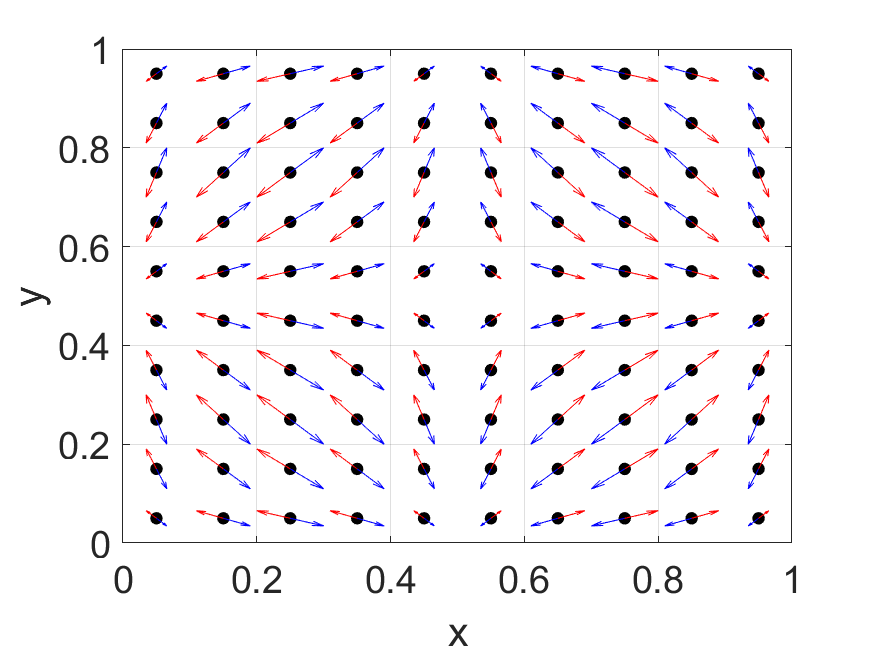}
\end{minipage}
\begin{minipage}[b]{0.49\textwidth}
	\includegraphics[width=\textwidth]{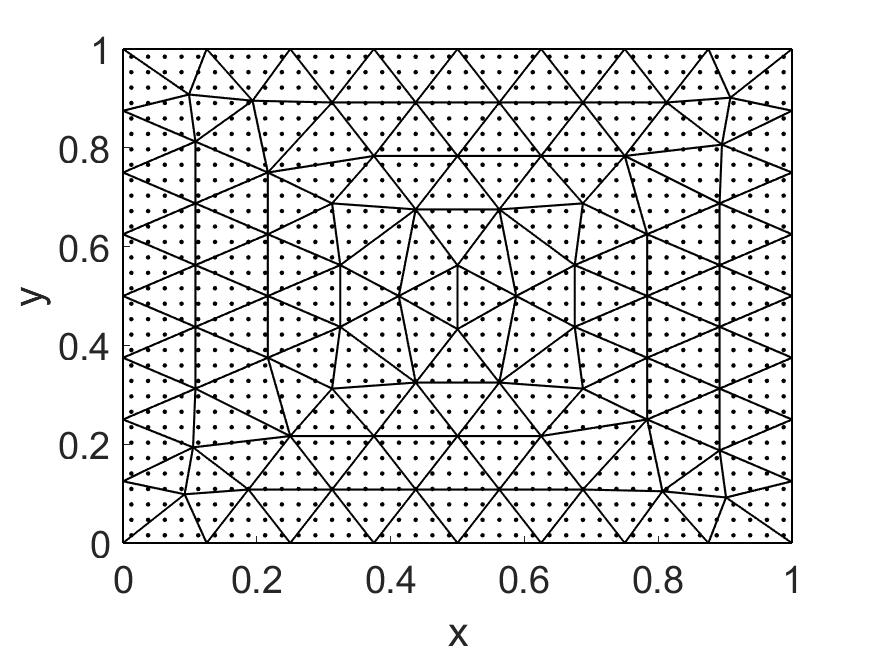}
\end{minipage}
\caption{The exact solution in which particles (marked with dots) move back and forth along the marked vectors (left) and the unstructured grid with the initial particle configuration (right).}
\label{fig:QuiverAndGrids}
\end{figure}

This benchmark was simulated with standard MPM and PS-MPM, using an unstructured triangular grid with material points initialised uniformly over the domain, as shown in Figure~\ref{fig:QuiverAndGrids}. A time step size of $\Delta t = 2.25 \cdot 10^{-4}$ s was chosen, which results in a Courant number of approximately $0.36$. For the adopted parameters, the imposed solution in Equation~\ref{eq:MMS_ux}-\ref{eq:MMS_uy} has period of $T=0.02$ s. 

\subsubsection{Grid-crossing error}

\begin{figure}[htb]
	\centering
	\begin{minipage}[b]{0.4\textwidth}
		\includegraphics[width=\textwidth]{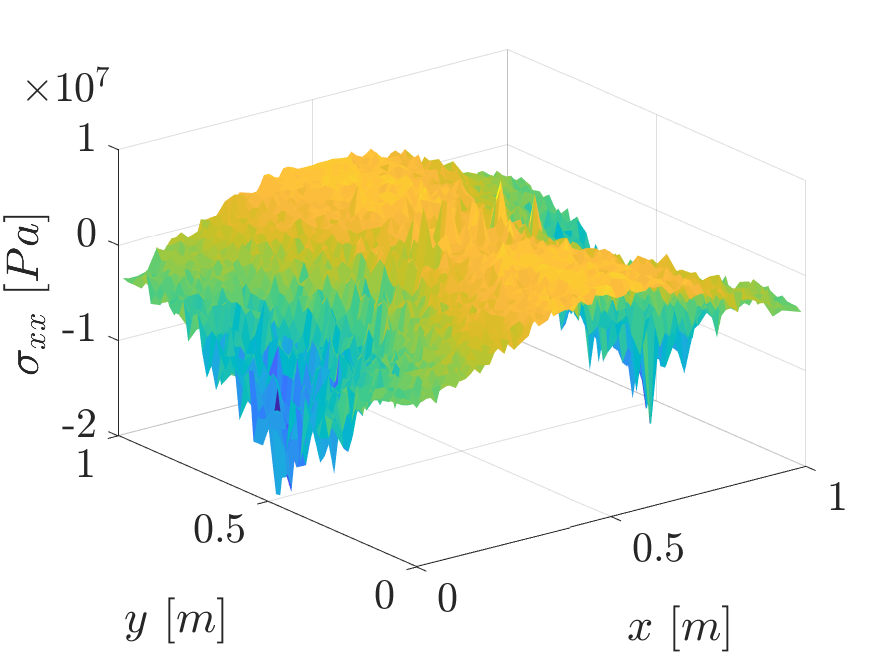}
		\centering{Standard MPM}
	    \label{fig:MMSStressField_ST}
	\end{minipage}
	\begin{minipage}[b]{0.4\textwidth}
		\includegraphics[width=\textwidth]{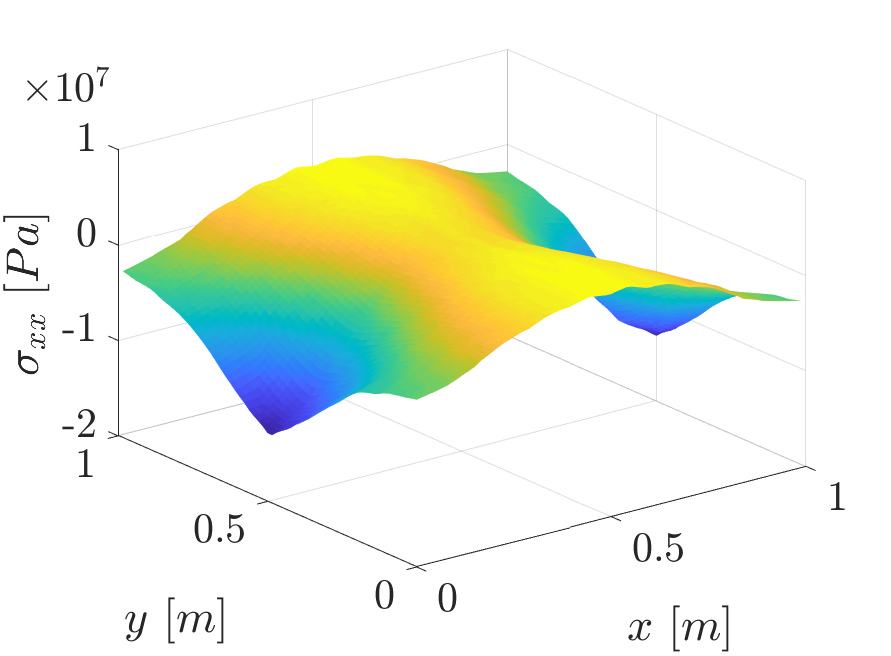}
		\centering{PS-MPM}
	    \label{fig:MMSStressField_PS}
	\end{minipage}
	\begin{minipage}[b]{0.4\textwidth}
		\includegraphics[width=\textwidth]{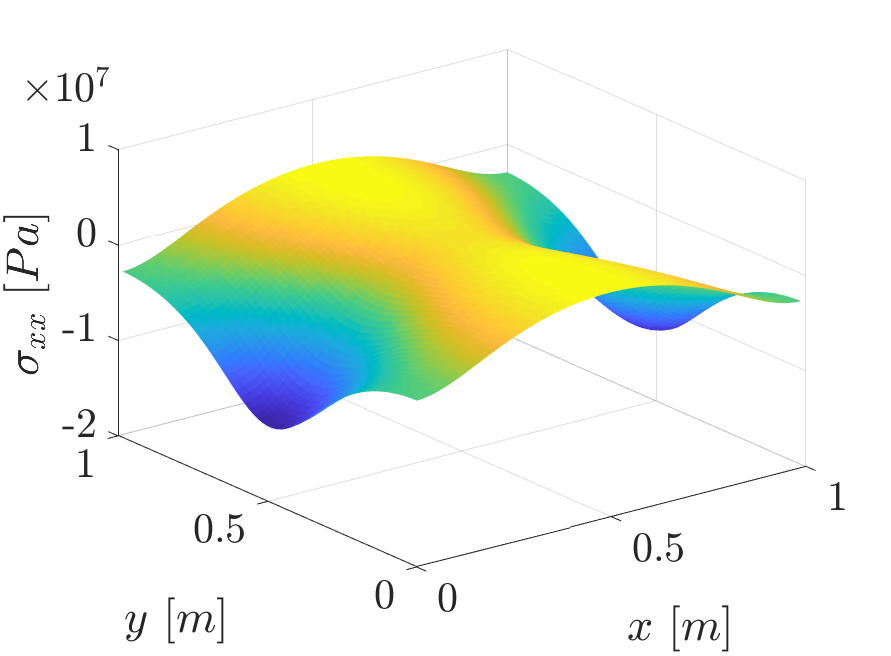}
		\centering{Exact}
	    \label{fig:MMSStressField_exact}
	\end{minipage}
	\caption{The interpolated particle stress in the $x$-direction at $t=0.016$s.}
	\label{fig:MMSStressField}
\end{figure}

\begin{figure}[tbh]
	\centering
    \begin{minipage}[b]{0.48\textwidth}
		\includegraphics[width=\textwidth]{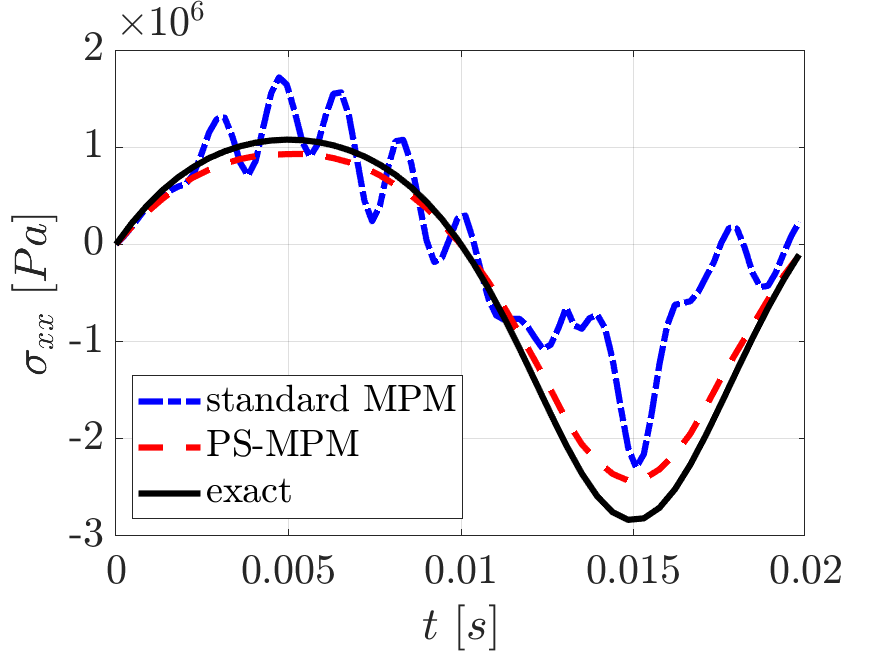}
	\end{minipage}
	\begin{minipage}[b]{0.48\textwidth}
		\includegraphics[width=\textwidth]{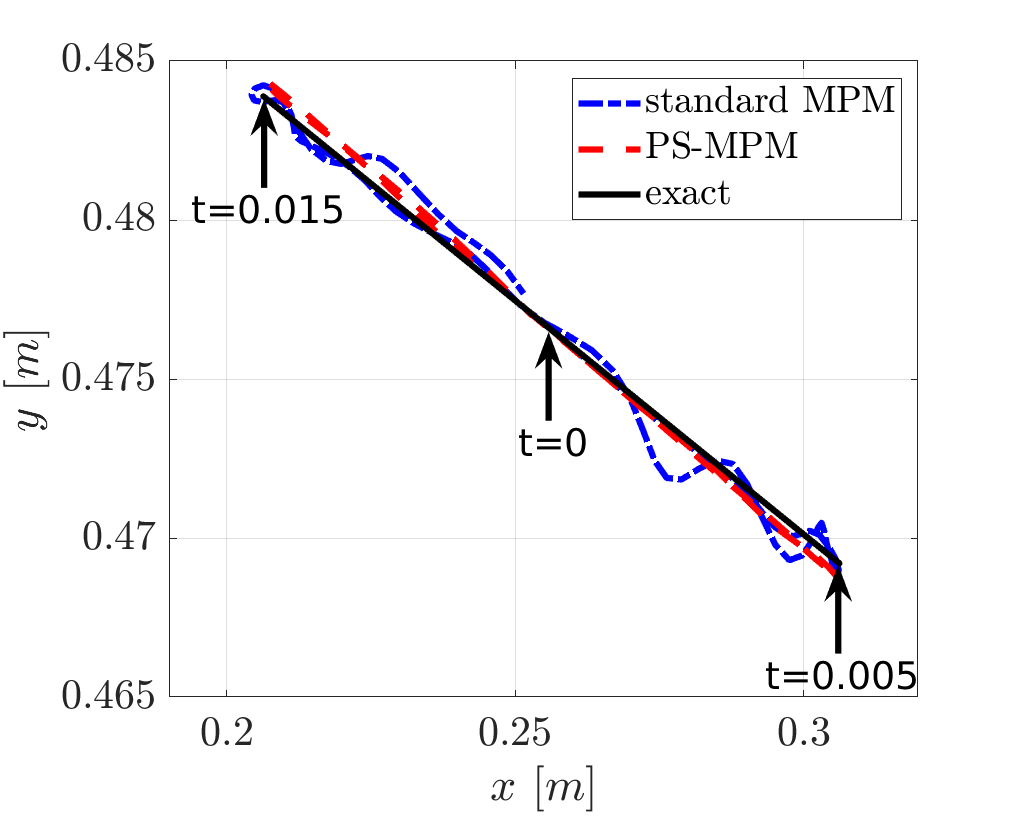}
	\end{minipage}
    \caption{The stress (left) and movement (right) of a single particle over time (in seconds) for both standard MPM and PS-MPM compared to the exact solution.}
	\label{fig:ParticleMovementAndStress}
\end{figure}

First, it is shown that the grid-crossing error typical for standard MPM does not occur in PS-MPM, by comparing the normal horizontal stress resulting from these methods (see Figure~\ref{fig:MMSStressField}). The configurations for PS-MPM and standard MPM contain the same number of particles ($5120$, 4 times as many as in Figure~\ref{fig:QuiverAndGrids}) and a comparable number of basis functions, $289$ for standard MPM and $243$ basis functions for PS-MPM. As expected, the stress field in standard MPM suffers severely from grid-crossing, whereas with PS-MPM a smooth stress field is observed.
The same conclusion can also be drawn when investigating individual particle stresses over time, as is shown in Figure~\ref{fig:ParticleMovementAndStress}. For this figure, the particle positioned at $(x^0,y^0)\approx(0.25$~m,$0.47$~m$)$ has been traced throughout the simulation. The figure depicts the stress and trajectory of the particle. For standard MPM, it is observed that the stress profile is highly oscillatory, resulting in a disturbed trajectory. PS-MPM shows no oscillations, but only a small offset, which was found to disappear upon further refinement of the grid.

\subsubsection{Spatial convergence}
PS-MPM with quadratic basis functions is expected to show third-order spatial convergence. To determine the spatial convergence, the time averaged root-mean-squared (RMS) error is used:
\begin{equation}\label{eq:RMS-Error}
    E(t) = \sqrt{\frac{\sum_{i=1}^{n_t}\sum_{p=1}^{n_p} |x_p(t_i)-\hat{x}_p(t_i)|^2}{n_p n_t}},
\end{equation}
where $n_t$ denotes the total number of time steps of the simulation.
Under the assumption that the spatial error is much larger than the error produced by the numerical integration and time-stepping scheme, the spatial convergence of standard MPM and PS-MPM is determined by varying the typical element length $h$. This length which was defined as the average edge length for standard MPM and the average sub-triangle edge length for PS-MPM.
It has been observed that the time-integration error is indeed sufficiently small, but the number of particles required for an adequately accurate integration increases rapidly as $h$ decreases. 
\begin{figure}[ht!]
	\centering
	\begin{tabular}{ c c}
	  \text{Standard MPM} & \text{PS-MPM}\\
	\begin{minipage}[b]{0.46\textwidth}
		\includegraphics[width=\textwidth]{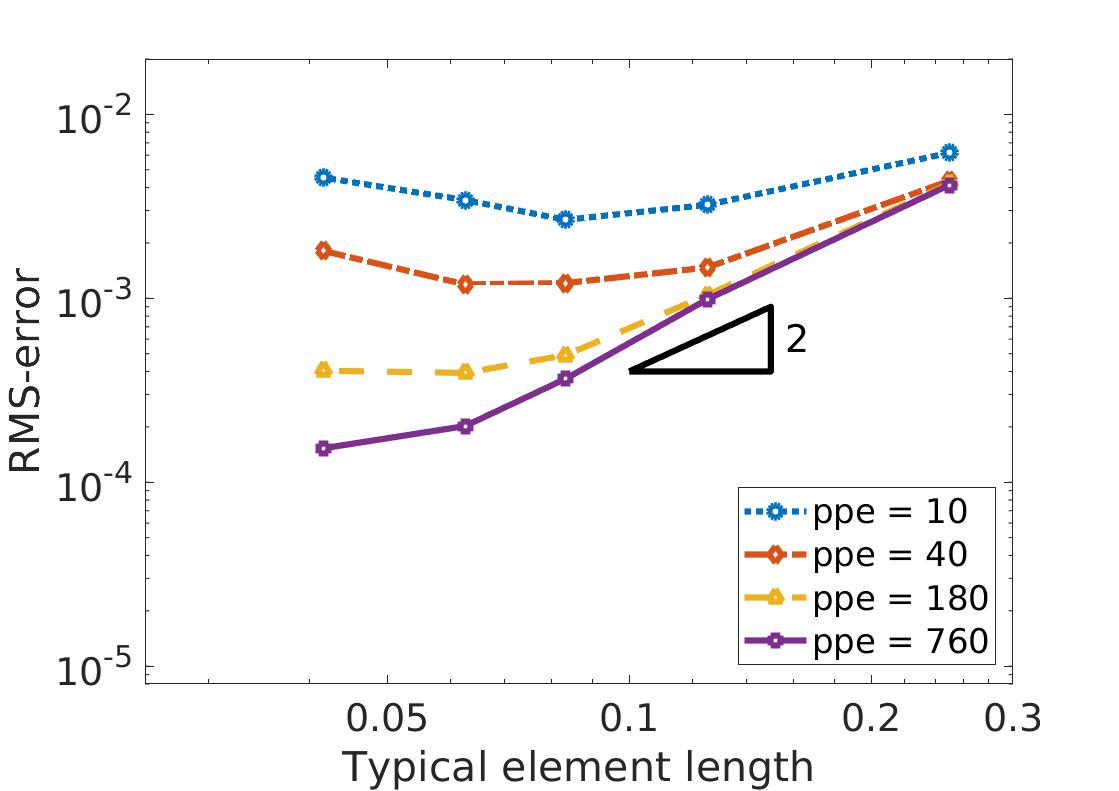}
	\end{minipage}
	&
	\begin{minipage}[b]{0.46\textwidth}
		\includegraphics[width=\textwidth]{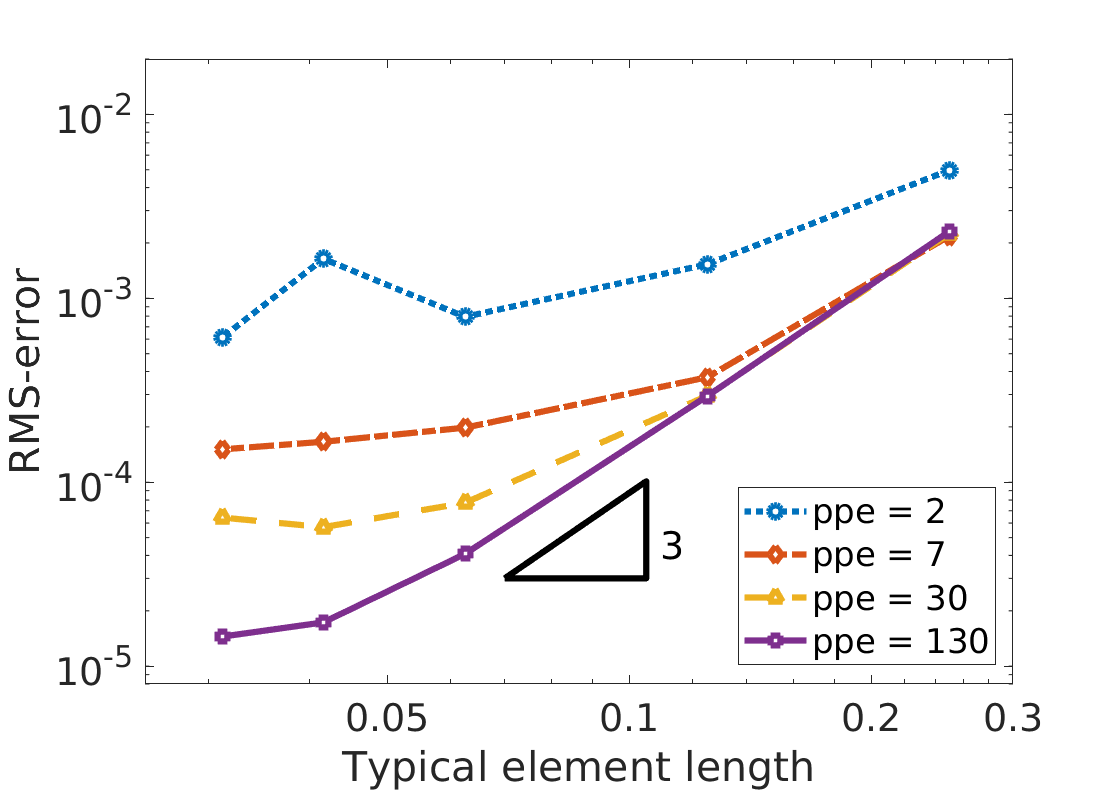}
	\end{minipage}\\
	\begin{minipage}[b]{0.46\textwidth}
		\includegraphics[width=\textwidth]{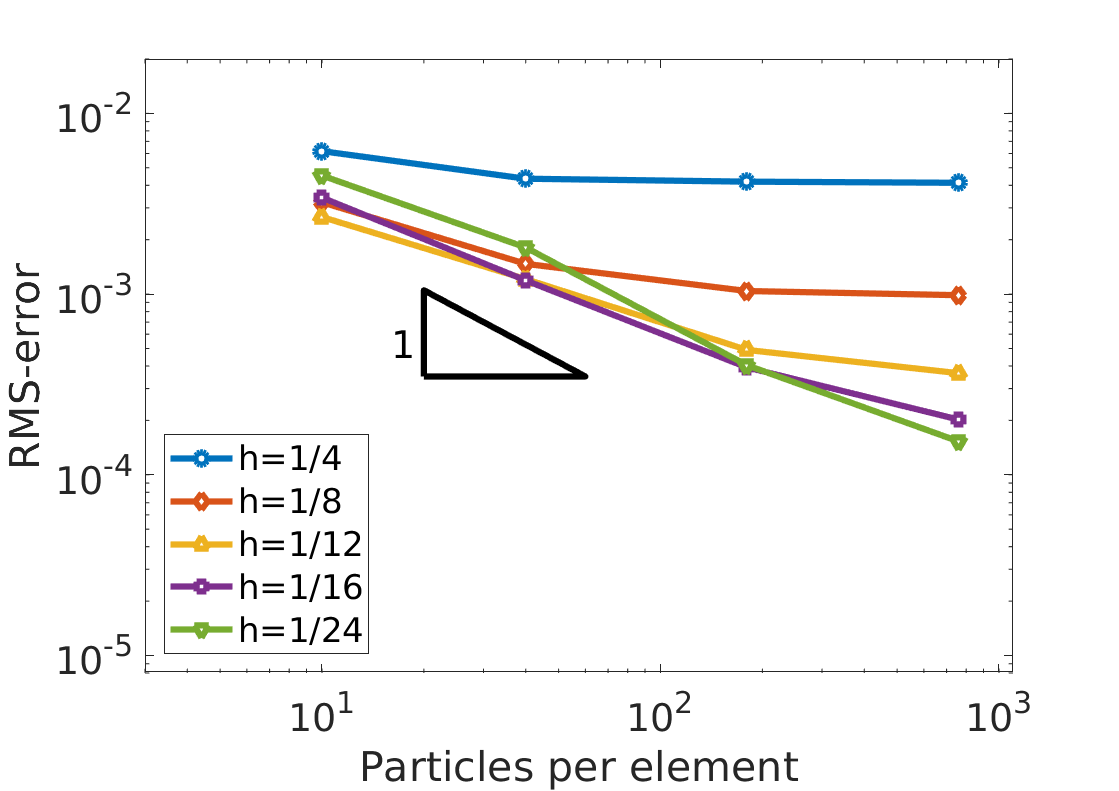}
	\end{minipage} 
	&
	\begin{minipage}[b]{0.46\textwidth}
		\includegraphics[width=\textwidth]{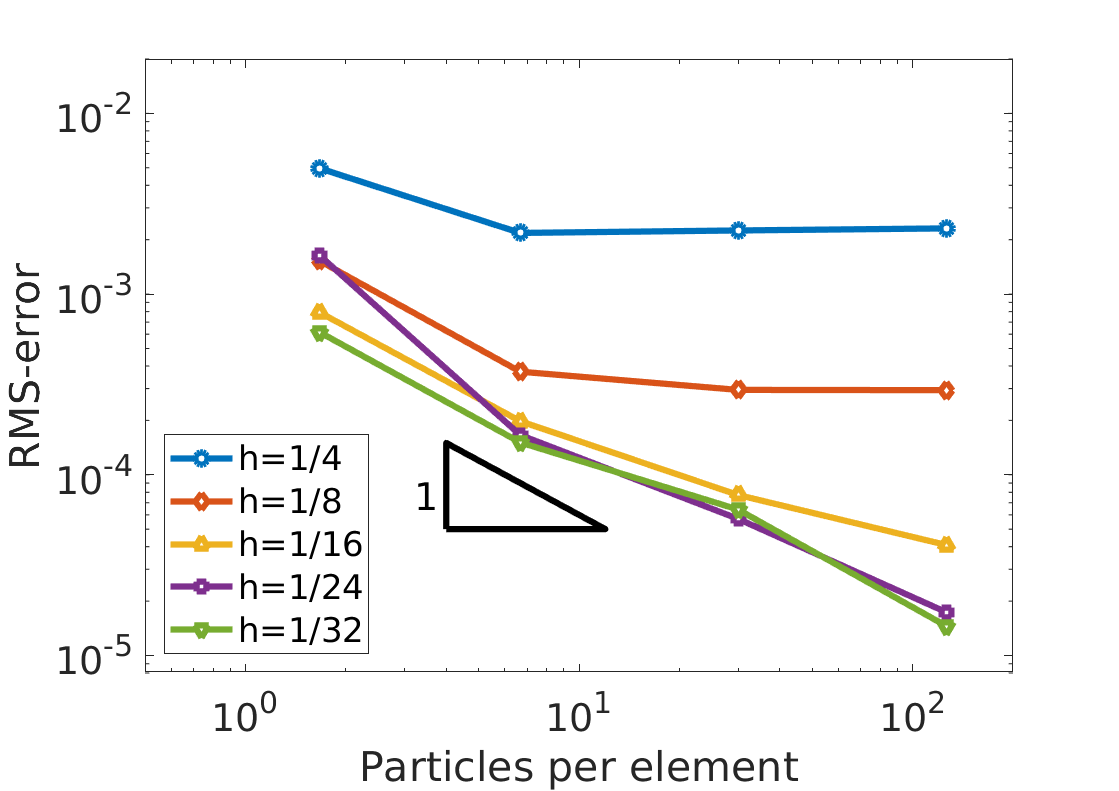}
	\end{minipage}
	\end{tabular}
    \caption{Spatial convergence (top) and convergence in the number of particles per element (ppe) (bottom) for both standard MPM (left) and PS-MPM (right).}
	\label{fig:MMS_Convergence}
\end{figure}

Figure~\ref{fig:MMS_Convergence} shows the spatial convergence of both standard MPM and PS-MPM for different numbers of particles per element on an unstructured grid. Provided that a sufficient number of particles is adopted, standard MPM shows second-order spatial convergence, whereas PS-MPM converges with third order. Furthermore, the RMS-error is smaller associated with PS-MPM for all considered configurations. To achieve the optimal convergence order for both standard MPM and PS-MPM, the number of integration particles required increases rapidly, as the integration error otherwise dominates the total error, as shown in Figure~\ref{fig:MMS_Convergence}. Note that the convergence rate in the number of particles for both standard MPM and PS-MPM is measured to be of first order. However, for a fixed typical element length $h$ and number of particles per element, the use of PS-MPM leads to a lower RMS-error, illustrating the higher accuracy per degree of freedom.  
The inaccurate integration in MPM due to the use of particles as integration points is known to limit the spatial convergence. Over the years, different measures have been proposed to decrease the quadrature error in MPM based on function reconstruction techniques like MLS \cite{gong2015improving}, cubic splines \cite{tielen2017high} and Taylor least squares \cite{wobbestaylor}. 
The combination of PS-MPM with function reconstruction techniques to obtain optimal convergence rates with a moderate number of particles is subject to future research to make PS-MPM even more suited for practical applications. 

\subsection{Column under self-weight}

In the previous examples, we considered PS-MPM adopting a consistent mass matrix. However, lumping of the mass matrix is common practice in many applications of MPM, as it speeds up the simulation and increases the numerical stability. In case elements become (almost) empty, the consistent mass matrix in Equation~\ref{eq:MassMatrixAssembled} has very small entries and becomes ill-conditioned, leading to stability issues, which is a well-known phenomena in MPM \cite{sulsky1994particle}. Using a lumped version of the mass matrix to solve for the acceleration and velocity fields, is known to overcome the ill-conditioning \cite{sulsky1994particle}. However, we have observed that lumping within PS-MPM also causes oscillations as a side-effect, which is not the case in standard MPM. 
\begin{figure}[htb]
    \centering
    \begin{minipage}[t]{0.22\textwidth}
		\begin{center}
			\def\svgwidth{1.0\columnwidth}
\begingroup%
  \makeatletter%
  \providecommand\color[2][]{%
    \errmessage{(Inkscape) Color is used for the text in Inkscape, but the package 'color.sty' is not loaded}%
    \renewcommand\color[2][]{}%
  }%
  \providecommand\transparent[1]{%
    \errmessage{(Inkscape) Transparency is used (non-zero) for the text in Inkscape, but the package 'transparent.sty' is not loaded}%
    \renewcommand\transparent[1]{}%
  }%
  \providecommand\rotatebox[2]{#2}%
  \newcommand*\fsize{\dimexpr\f@size pt\relax}%
  \newcommand*\lineheight[1]{\fontsize{\fsize}{#1\fsize}\selectfont}%
  \ifx\svgwidth\undefined%
    \setlength{\unitlength}{126.72306463bp}%
    \ifx\svgscale\undefined%
      \relax%
    \else%
      \setlength{\unitlength}{\unitlength * \real{\svgscale}}%
    \fi%
  \else%
    \setlength{\unitlength}{\svgwidth}%
  \fi%
  \global\let\svgwidth\undefined%
  \global\let\svgscale\undefined%
  \makeatother%
  \begin{picture}(1,2.86377581)%
    \lineheight{1}%
    \setlength\tabcolsep{0pt}%
    \put(0,0){\includegraphics[width=\unitlength,page=1]{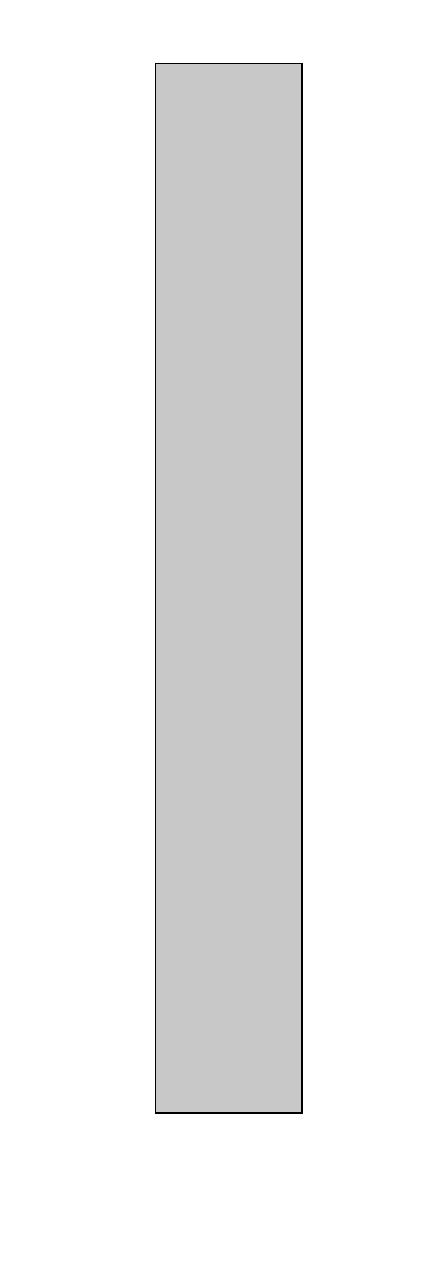}}%
    \put(0.46053926,1.50655187){\color[rgb]{0,0,0}\makebox(0,0)[lt]{\lineheight{0}\smash{\begin{tabular}[t]{l}$g$\end{tabular}}}}%
    \put(0.20448183,0.20261283){\color[rgb]{0,0,0}\makebox(0,0)[lt]{\lineheight{0}\smash{\begin{tabular}[t]{l}$0$\end{tabular}}}}%
    \put(0.71720343,0.20364329){\color[rgb]{0,0,0}\makebox(0,0)[lt]{\lineheight{0}\smash{\begin{tabular}[t]{l}$W$\end{tabular}}}}%
    \put(0.17790262,2.74969222){\color[rgb]{0,0,0}\makebox(0,0)[lt]{\lineheight{0}\smash{\begin{tabular}[t]{l}$H$\end{tabular}}}}%
    \put(0,0){\includegraphics[width=\unitlength,page=2]{SoilColumn2D.pdf}}%
    \put(0.42999574,0.01324421){\color[rgb]{0,0,0}\makebox(0,0)[lt]{\lineheight{0}\smash{\begin{tabular}[t]{l}$x$\end{tabular}}}}%
    \put(-0.00748012,1.48878313){\color[rgb]{0,0,0}\makebox(0,0)[lt]{\lineheight{0}\smash{\begin{tabular}[t]{l}$y$\end{tabular}}}}%
    \put(0,0){\includegraphics[width=\unitlength,page=3]{SoilColumn2D.pdf}}%
  \end{picture}%
\endgroup%

		\end{center}
	\end{minipage}
	\hspace{0.03\linewidth}
	\begin{minipage}[t]{0.35\linewidth}
		\centering
		\includegraphics[width=0.4\textwidth]{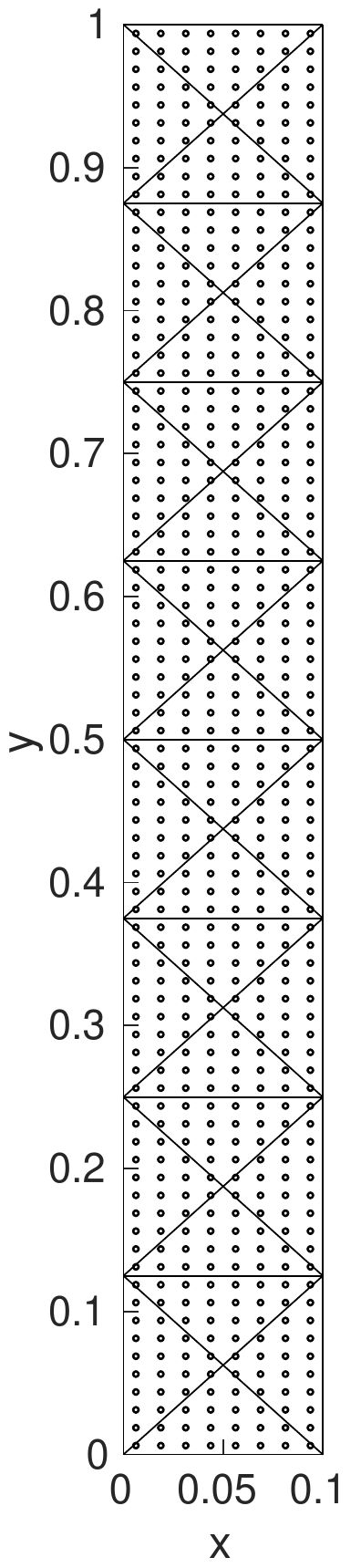}
		\includegraphics[width=0.4\textwidth]{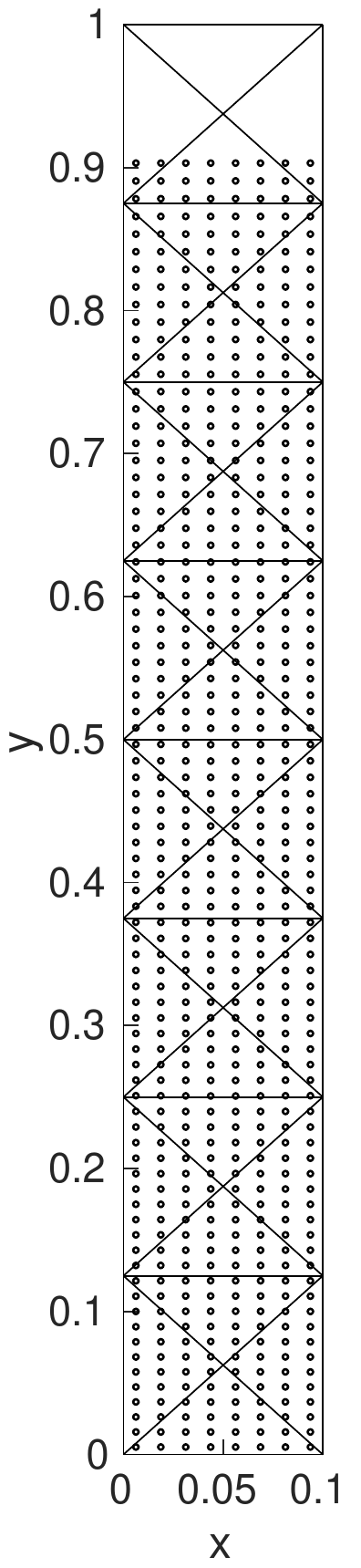}
	\end{minipage}
	\caption{The soil column and a discretisation with the expected material point positions at $t=0$ and $t=0.186$.}
	\label{fig:GridSoilColumn}
\end{figure}

To illustrate this phenomenon, a soil-column benchmark under self-weight is considered, as shown in Figure~\ref{fig:GridSoilColumn}. The soil column in this exemplary benchmark is modelled as a linear-elastic material. It is fixed in all directions at the bottom, and completely free at the top. At the left and right boundary, we impose a free-slip boundary condition for movement in the $y$-direction, but the displacement in the $x$-direction is fixed to be zero. The column is compressed by gravity due to its self-weight. 
The soil column is modelled with density $\rho=1\cdot 10^3\,\si{kg/m^3}$, Young's modulus $E=1\cdot 10^5\,\si{Pa}$, Poisson ratio $\nu=0$, gravitational acceleration $g=-9.81\, \si{m/s^2}$, height $H=1\,\si{m}$ and width $W=0.1\,\si{m}$. The maximum strain when adopting these parameters is approximately $18\%$.
\begin{figure}[h!tb]
	\centering
	\begin{tabular}{c c c}
	& 
	\text{PS-MPM, fully lumped} 
	& 
	\text{\hspace{0.5cm}PS-MPM, partially lumped}\\
	\rotatebox{90}{\hspace{0.3cm} Displacement middle particle} & 
	\begin{minipage}[b]{0.42\textwidth}
		\includegraphics[width=1.0\textwidth]{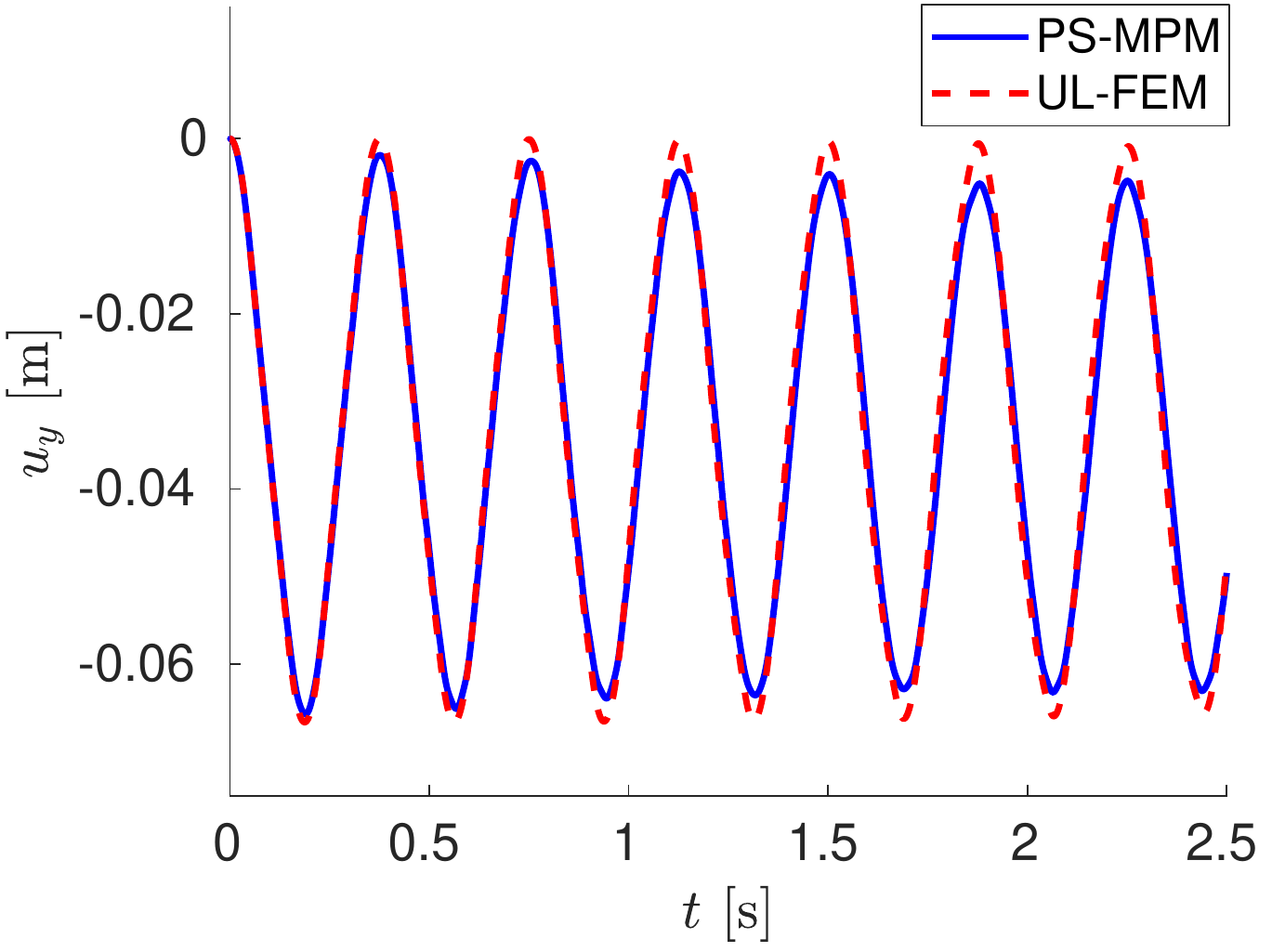}
	\end{minipage}
	&
	\begin{minipage}[b]{0.42\textwidth}
		\includegraphics[width=1.0\textwidth]{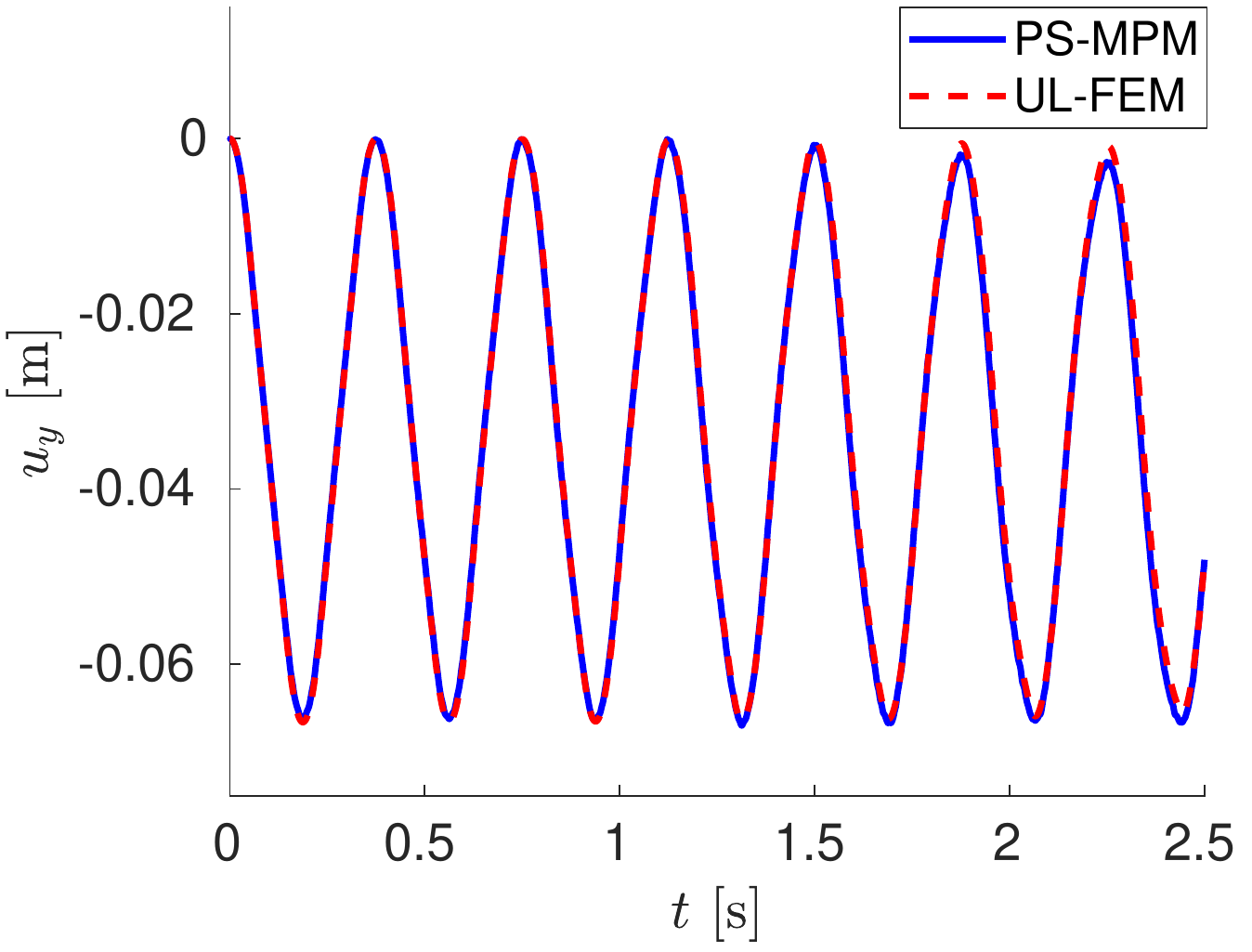}
	\end{minipage}\\
	\rotatebox{90}{\hspace{0.7cm}Velocity middle particle} &
	\begin{minipage}[b]{0.42\textwidth}
		\includegraphics[width=1.0\textwidth]{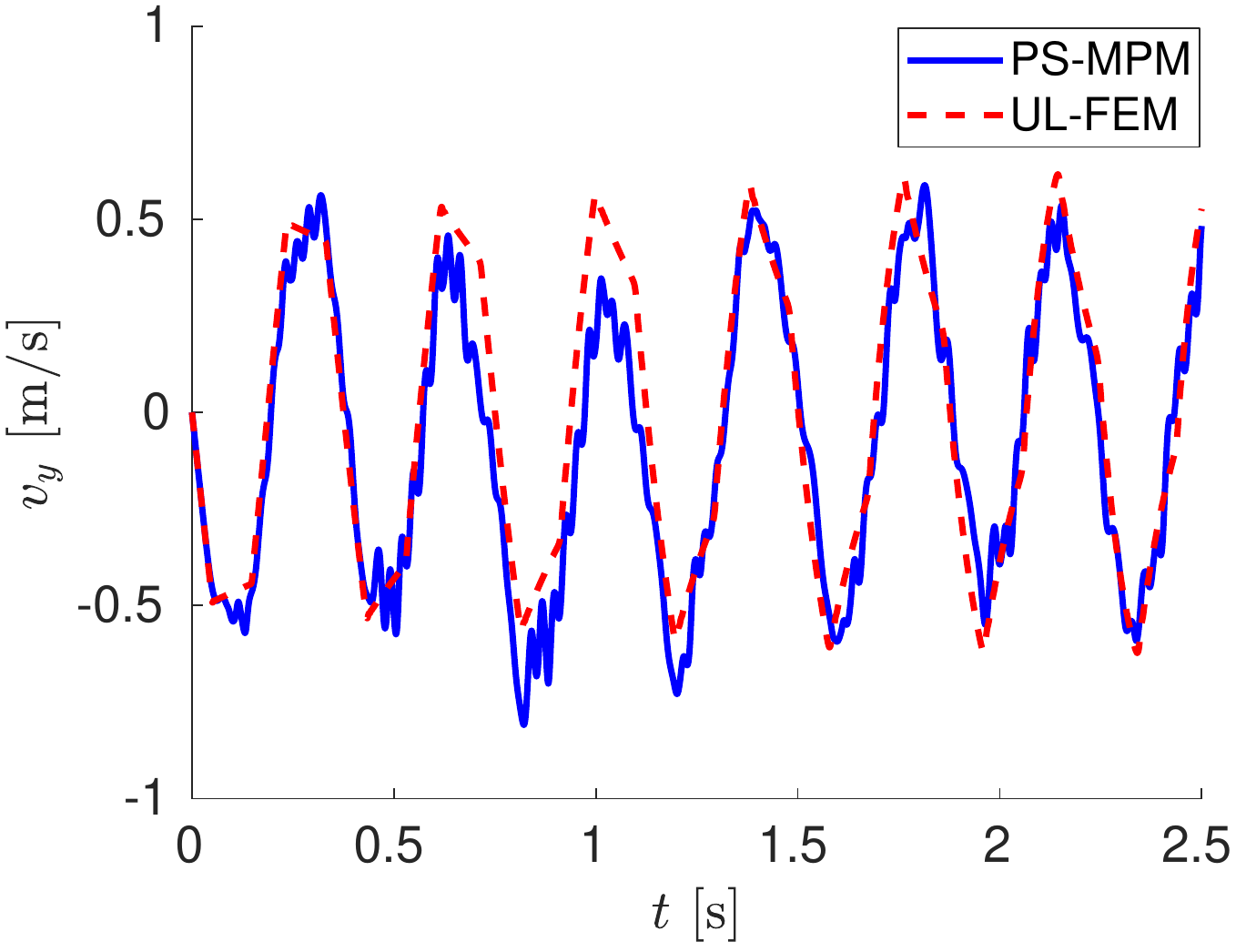}
	\end{minipage} 
	&
	\begin{minipage}[b]{0.42\textwidth}
	\includegraphics[width=1.0\textwidth]{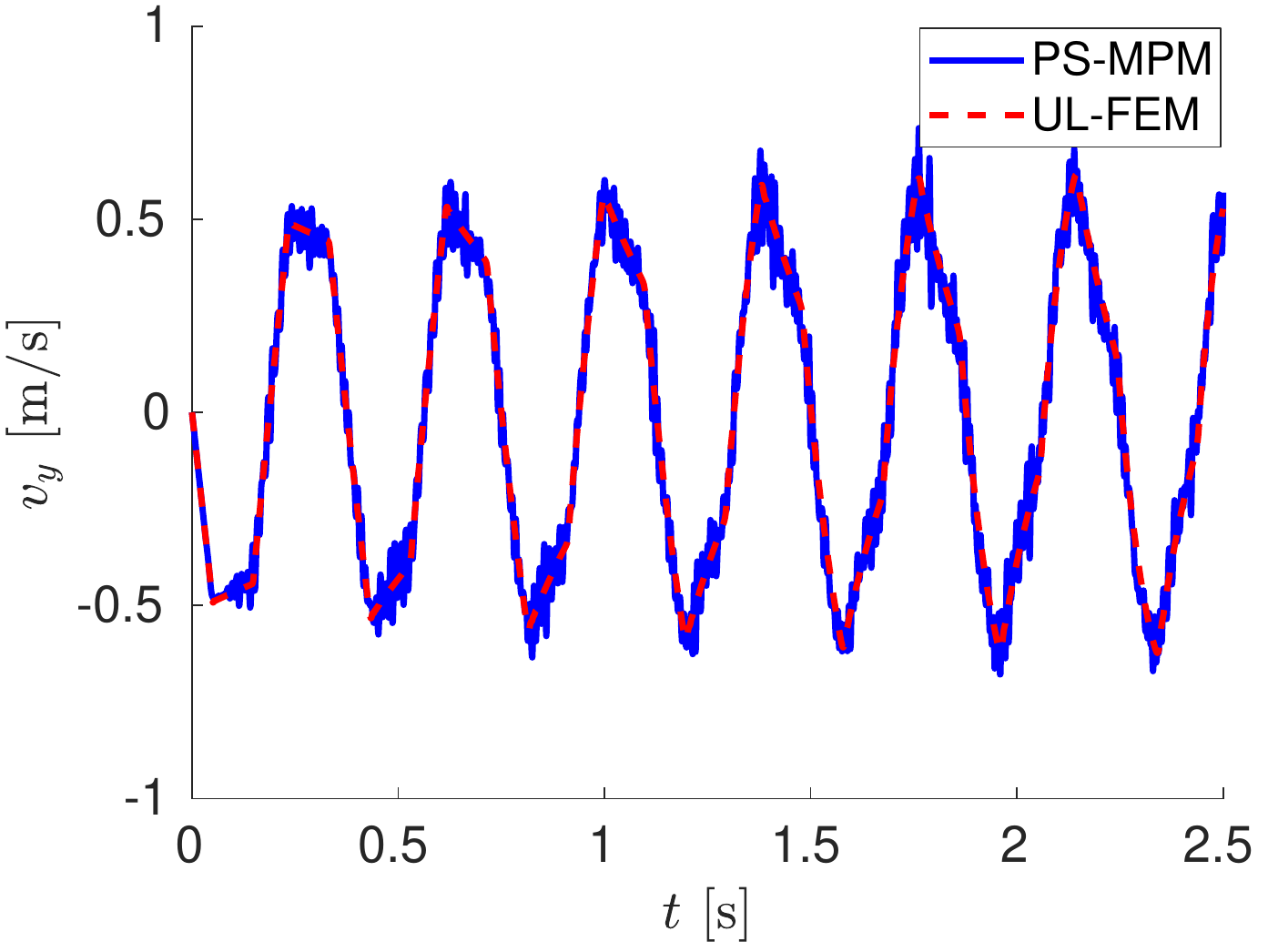}
	\end{minipage}\\
	\rotatebox{90}{\hspace{0.7cm}Stress profile at $t=2.5$s} &
	\begin{minipage}[b]{0.42\textwidth}
	\includegraphics[width=1.0\textwidth]{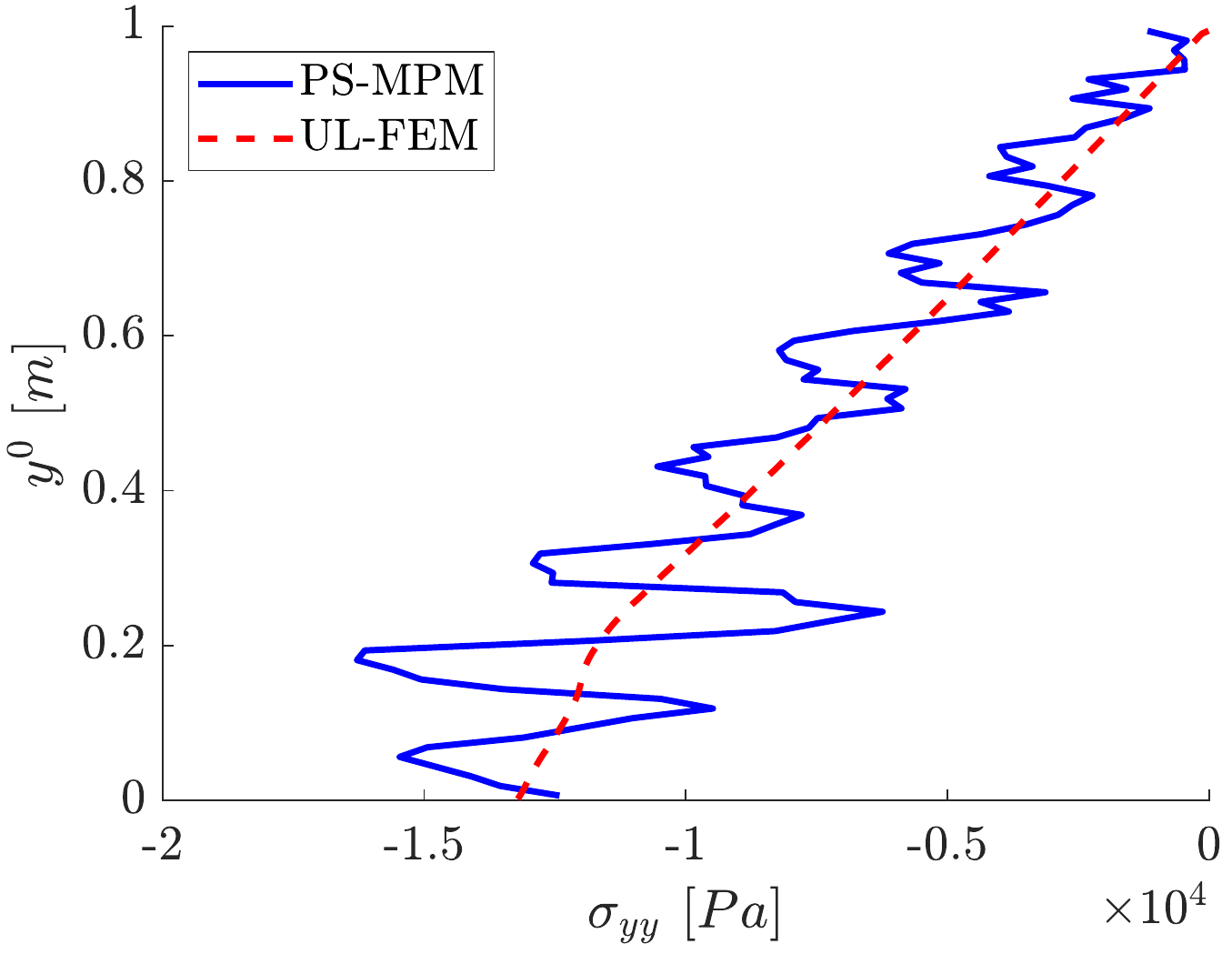}
	\end{minipage} 
	&
	\begin{minipage}[b]{0.42\textwidth}
	\includegraphics[width=1.0\textwidth]{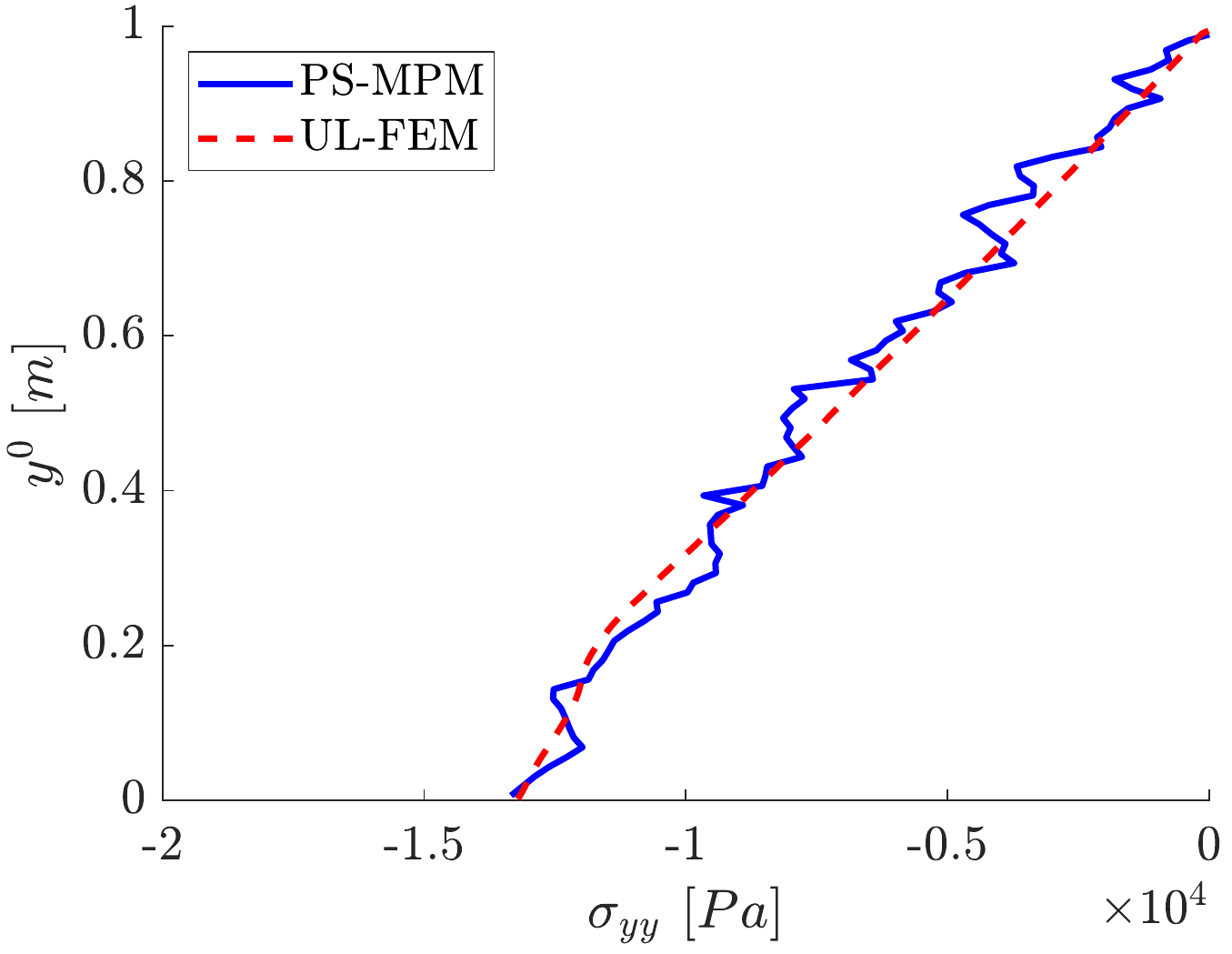}	\end{minipage}
	\end{tabular}
    \caption{Displacement and velocity of the middle particle and stress profile at $t=2.5$s with fully lumped (left) and partially lumped (right) PS-MPM.}
	\label{fig:SC_LD_lumped}
\end{figure}
Figure~\ref{fig:GridSoilColumn} illustrates the grid considered for this benchmark as well as the particle positions when the deformation is maximal. 

Due to the fact that part of the grid becomes empty, the use of a consistent mass matrix leads to stability issues with unlumped PS-MPM, which makes a PS-MPM simulation impossible. Therefore, a lumped mass matrix was adopted for this benchmark to restore stability. The use of a lumped mass matrix within PS-MPM indeed restores the stability, but spatial oscillations occur in the stress profile, degrading the solution of the displacement and velocity as well, as shown in Figure~\ref{fig:SC_LD_lumped} (left column). 

\subsubsection{Partial lumping to mitigate spatial oscillations}
A solution to mitigate the oscillations when lumping the mass matrix in PS-MPM, is the application of partial lumping. Instead of applying full lumping on the mass matrix, partial lumping only lumps those rows responsible for the ill-conditioning. These rows typically correspond to the basis functions in the part of the grid where very few particles are left. For this benchmark, all rows corresponding to a basis function with at least one empty main-triangle in its molecule are lumped. 
Additional empty elements are added at the top of the column, to ensure that the top-most basis functions are always lumped.     
Figure~\ref{fig:SC_LD_lumped} (right column) shows the displacement and velocity of a particle over time, as well as the stress profile at $t=2.5$~s obtained with partial lumped PS-MPM. Compared to the results obtained with lumped PS-MPM, shown in Figure~\ref{fig:SC_LD_lumped} (left), the stress profile, the velocity and displacement over time significantly improve. 

Hence, the use of partial lumping within PS-MPM combines the advantages of adopting a consistent and lumped mass matrix, while minimising their drawbacks, improving the overall performance of the method. Future research will focus on techniques to further mitigate or fully overcome the oscillations caused by lumping, in particular when PS-MPM is applied for practical problems. 

\section{Conclusion}\label{sec:Conclusion}
	
In this paper, we presented an alternative for B-spline MPM suited for unstructured triangulations using piecewise quadratic $C^1$-continuous Powell-Sabin spline basis functions. The method combines the benefits of smooth, higher-order basis functions with the geometric flexibility of an unstructured triangular grid. 
PS-MPM yields a mathematically sound approach to eliminate grid-crossing errors, due to the smooth gradients of the basis functions. As a first validation, a vibrating bar was considered, for which the PS-MPM solution yields accurate results on a relatively coarse grid. Numerical simulations, obtained for a vibrating plate undergoing axis-aligned displacement, have shown higher-order convergence for the particle stresses and displacements. 

The use of a lumped mass matrix to increase stability within PS-MPM, leads to spatial oscillations in the stress profile. A partial lumping strategy was proposed, to combine the advantages of adopting a consistent and lumped mass matrix, which successfully mitigates this issue. Investigation of alternative unstructured spline technologies, in particular, refinable $C^1$ splines on irregular quadrilateral grids~\cite{Nguyen2016} is underway.

\section*{Conflict of interest}

Pascal de Koster has received funding support from Dutch research institute Deltares.

\bibliographystyle{spmpsci}      
\bibliography{References}   

\end{document}